\documentclass[11pt]{amsart}
\usepackage{graphicx}
\def\t{\widetilde}
\def\bb{\Bbb}
\def\a{\alpha}        
\def\b{\beta}        
\def\c{\gamma}           
\def\d{\delta}           
\def\e{\varepsilon}      \def\x{\chi}
\def\f{\zeta}            \def\y{\psi}
\def\g{\eta}             \def\z{\zeta}  
\def\h{\theta}    
\def\Ga{\Gamma}  
\def\De{\Delta}  
\def\Th{\Theta}  
 
\def\Xi{\Xi}     
\def\Pi{\Pi}     
\def\A{\mathcal{A}}  \def\F{\mathcal{F}}
\def\B{\mathcal{B}}  \def\G{\mathcal{G}}
\def\C{\mathcal{C}}  \def\H{\mathcal{H}}
\def\D{\mathcal{D}}  \def\I{\mathcal{I}}
\def\E{\mathcal{E}}  \def\J{\mathcal{J}}
\newtheorem{Thm}{Theorem}[section]
\newtheorem{Cor}[Thm]{Corollary}

\newtheorem{Prop}[Thm]{Proposition}
\newtheorem{Lem}[Thm]{Lemma}
\newtheorem{Claim}[Thm]{Claim}
   \title{A Criterion for Almost Alternating Links to be
          Non-splittable}
  \author{Tatsuya TSUKAMOTO}
 \address{Department of Mathematics, School of Science and
          Engineering, Waseda University, 3-4-1 Okubo, 
          Shinjuku-ku, Tokyo 169-8555, Japan}
\curraddr{Department of Mathematics, the George Washington
          University, Funger Hall 428, 2201 G St. N.W. 
          Washington, D.C. 20052, USA}
   \email{tatsuya@research.circ.gwu.edu}
  \thanks{The author is partially supported by JSPS Research
          Fellowships for Young Scientists.}
  \thanks{The paper is dedicated to Professor Shin'ichi Suzuki 
          on his sixtieth birthday.}
\date{}
\begin{document}
\maketitle
%%%%%%%%%%%%%%%%
%
%%%%%%%%%%%%%%%%%%%%%%%%%%%%%%%%%%%%%%%%%%%%%% INTRODUCTION %%%%%
%
%%%%%%%%%%%%%%%%%%%%%%
\section{Introduction}
%%%%%%%%%%%%%%%%%%%%%%
%
%%%%%%%%%%%%%%%%%%%%%%%%%%%%%%%%%%%
The notion of almost alternating links was introduced by 
C. Adams et al (\cite{a1}). Here we give a sufficient condition 
for an almost alternating link diagram to represent a
non-splittable link. This solves a question asked in \cite{a1}. 
A partial solution for special almost alternating links has been 
obtained by M. Hirasawa (\cite{h}). 
%%%%%%%%%%%%%%%%%%%%%%%%%%%%%%%%%%%
%
%%%%%%%%%%%%%%%%%%%%%%%%%%%%%%%%%%%
As its applications, Theorem \ref{ThmLinks} gives us a way to 
see if a given almost alternating link diagram represents a 
splittable link without increasing numbers of crossings of 
diagrams in the process. Moreover, we show that almost 
alternating links with more than two components are non-trivial.
In Section \ref{SecResults}, we state them in detail.
%%%%%%%%%%%%%%%%%%%%%%%%%%%%%%%%%%%
%
%%%%%%%%%%%%%%%%%%%%%%%%%%%%%%%%%%%
To show our theorem, we basically use a technique invented 
by W. Menasco (see \cite{m, mt}). We review it in Section 
\ref{SecStdPosition}. However, we also apply `` charge and 
discharge method" to our graph-theoretic argument, which is 
used to prove the four color theorem in \cite{ah}. 
%%%%%%%%%%%%%%%%%%%%%%%%%%%%%%%%%%%
%
%%%%%%%%%%%%%%%%
\par\vspace{5mm}
%%%%%%%%%%%%%%%%
%
%%%%%%%%%%%%%%%%%%%%%%%%%%%%%%%%%%%%%%%%%%%%%%%%%%% RESULTS %%%%%
%
%%%%%%%%%%%%%%%%%%%%%
%
%%%%%%%%%%%%%%%%%%%%%%%%%%%%%%%%%%%%%%%%%%%%%%%%%%
\section{The main theorem and its applications}\label{SecResults}
%%%%%%%%%%%%%%%%%%%%%%%%%%%%%%%%%%%%%%%%%%%%%%%%%%
%
%%%%%%%%%%%%%%%%%%%%%%%%%%%%%%%%%%%
Menasco has shown that an alternating link diagram can represent
a splittable link only in a trivial way.
%%%%%%%%%%%%%%%%%%%%%%%%%%%%%%%%%%%
%
%%%%%%%%%%%%%%%%%%%%%%%%%%%%%%%%%%%%%%%%% THEOREM ThmALinks %%%%%
\begin{Thm}\label{ThmALinks}$($\cite{m}$)$ If a link $L$ has a 
connected alternating diagram, then $L$ is non-splittable.  
\end{Thm}                                                       
%%%%%%%%%%%%%%%%%%%%%%%%%%%%%%%%%%%%%%%%%%%%%%%%%%%%%%%%%%%%%%%%%
%
%%%%%%%%%%%%%%%%%%%%%%%%%%%%%%%%%%%
We say a link diagram $\t{L}$ on $S^2$ is {\it almost alternating} 
if one crossing change makes $\t{L}$ alternating. A link $L$ is 
{\it almost alternating} if $L$ is not alternating and $L$ has an 
almost alternating diagram. We call a crossing of an almost 
alternating diagram a {\it dealternator} if the crossing change at 
the crossing makes the diagram alternating. An almost alternating 
diagram may have more than one dealternator. However, we can 
uniquely decide a connected almost alternating diagram if the 
diagram has more than one dealternater (Proposition 
\ref{PropTwoDealt}). Since the statement of Proposition 
\ref{PropTwoDealt} does not contradict the statement 
of Theorem \ref{ThmLinks}, we may assume that our almost 
alternating diagram has exactly one dealternator from now on.
%%%%%%%%%%%%%%%%%%%%%%%%%%%%%%%%%%%
%
%%%%%%%%%%%%%%%%%%%%%%%%%%%%%%%%%% PROPOSITION PropTwoDealt %%%%%
\begin{Prop}\label{PropTwoDealt}If a connected almost alternating 
diagram $\t{L}$ has more than one dealternator, then $\t{L}$ is a 
diagram obtained from a Hopf link diagram with two crossings by 
changing one of the crossings.           \end{Prop} 
%%%%%%%%%%%%%%%%%%%%%%%%%%%%%%%%%%%%%%%%%%% END PROPOSITION %%%%%
%
%%%%%%%%%%%%%%%%%%%%%%%%%%%%%%%%%%%%%%%%%%%%%%%%%%%%% PROOF %%%%%
\begin{proof}                                                   
Assume that $\t{L}$ has more than one dealternator.             
Let $\a$ be one of the dealternators. Then, $\a$ is adjacent to 
other four crossings (some of them may be the same).     
Let $\b$ be another dealternator. Since the crossing change at  
$\b$ makes $\t{L}$ alternating, each of those four crossings    
must coinsides $\b$. Then, $\t{L}$ is a diagram obtained        
from a Hopf link diagram with two crossings by changing one of  
the crossings.                                                  
\end{proof}                                                     
%%%%%%%%%%%%%%%%%%%%%%%%%%%%%%%%%%%%%%%%%%%%%%%%% END PROOF %%%%%
%
%%%%%%%%%%%%%%%%
\par\vspace{3mm}
%%%%%%%%%%%%%%%%
%
%%%%%%%%%%%%%%%%%%%%%%%%%%%%%%%%%%%
A diagram $\t{L}$ on $S^2$ is {\it prime} if $\t{L}$ is connected, 
$\t{L}$ has at least one crossing, and there does not exist a 
simple closed curve on $S^2$ meeting $\t{L}$ transversely in just 
two points belonging to different arcs of $\t{L}$. If our almost 
alternating diagram is non-prime, then it is a connected sum of a 
prime almost alternating diagram and alternating diagrams. 
Therefore, we may restrict our interest to prime diagrams, since 
we know that a connected alternating diagram represents 
a non-splittable link (Theorem \ref{ThmALinks}).
%%%%%%%%%%%%%%%%%%%%%%%%%%%%%%%%%%%
%
%%%%%%%%%%%%%%%%
\par\vspace{3mm}
%%%%%%%%%%%%%%%%
%
%%%%%%%%%%%%%%%%%%%%%%%%%%%%%%%%%%%
In this paper, we say a diagram is {\it reduced} if it is not 
diagram I or diagram II in Figure \ref{fig:NonRedDiag}, where we 
allow both two types of crossings which give almost alternating 
diagrams with the dealternator at the marked double point. Then our 
main theorem is the following.
%%%%%%%%%%%%%%%%%%%%%%%%%%%%%%%%%%%
%
%%%%%%%%%%%%%%%%%%%%%%%%%%%%%%%%%%%%%%%%%% THEOREM ThmLinks %%%%%
\begin{Thm}\label{ThmLinks} If a link $L$ has a connected, prime,
reduced almost alternating diagram, then $L$ is non-splittable.
\end{Thm}
%%%%%%%%%%%%%%%%%%%%%%%%%%%%%%%%%%%%%%%%%%%%%%% END THEOREM %%%%%
%
%%%%%%%%%%%%%%%%%%%%%%%%%%%%%%%%%%%
\begin{figure}[htbp]\caption{Non reduced diagrams}
\label{fig:NonRedDiag}\end{figure}\par\vspace{3mm}
%%%%%%%%%%%%%%%%%%%%%%%%%%%%%%%%%%%
%
%%%%%%%%%%%%%%%%%%%%%%%%%%%%%%%%%%%
Theorem \ref{ThmLinks} gives us a way to see if a given almost 
alternating diagram $\t{L}$ represents a splittable link or not; 
here we may assume that $\t{L}$ is connected and prime. If 
$\t{L}$ is reduced, then it represents a non-splittable link 
from Theorem \ref{ThmLinks}. If $\t{L}$ is not reduced, then it 
is diagram I or diagram II in Figure \ref{fig:NonRedDiag}. If 
it is diagram I, then we obtain the diagram with no crossings of 
the trivial link with two components or an alternating diagram 
from $\t{L}$ by applying reducing move I to $\t{L}$ (see Figure 
\ref{fig:RedMoves}). Then, we can see if $\t{L}$ represents a 
splittable link or not from Theorem \ref{ThmALinks}. If it is 
diagram II, then we obtain another connected, prime almost 
alternating diagram $\t{L'}$ from $\t{L}$ with less crossings 
than those of $\t{L}$ by applying reducing move II to $\t{L}$. 
Then, we can see if $\t{L}$ represents a splittable link or 
not by continuing this process as far as we have diagram II.
%%%%%%%%%%%%%%%%%%%%%%%%%%%%%%%%%%%
%
%%%%%%%%%%%%%%%%%%%%%%%%%%%%%%%%%%%
\par\vspace{3mm}
\begin{figure}[htbp]\caption{Reducing moves}
\label{fig:RedMoves}\end{figure}\par\vspace{3mm}
%%%%%%%%%%%%%%%%%%%%%%%%%%%%%%%%%%%
%
%%%%%%%%%%%%%%%%%%%%%%%%%%%%%%%%%%%
As a corollary of Theorem \ref{ThmLinks}, we immediately have 
the following.
%%%%%%%%%%%%%%%%%%%%%%%%%%%%%%%%%%%
%
%%%%%%%%%%%%%%%%%%%%%%%%%%%%%%%%%%%%%% COROLLARY CorTrivial %%%%%
\begin{Cor}\label{CorTrivial} If a link $L$ has a connected, 
prime, reduced almost alternating diagram, then $L$ is 
non-trivial.
\end{Cor}
%%%%%%%%%%%%%%%%%%%%%%%%%%%%%%%%%%%%%%%%%%%%% END COROLLARY %%%%%
%
%%%%%%%%%%%%%%%%%%%%%%%%%%%%%%%%%%%
Moreover, we obtain the following if we restrict our interest 
to links with more than two components.
%%%%%%%%%%%%%%%%%%%%%%%%%%%%%%%%%%%
%
%%%%%%%%%%%%%%%%%%%%%%%%%%%%%%%%%%%% COROLLARY CorThreeComp %%%%%
\begin{Cor}\label{CorThreeComp} If a link $L$ with more than two 
components has a connected almost alternating diagram, then $L$ 
is non-trivial.       \end{Cor}
%%%%%%%%%%%%%%%%%%%%%%%%%%%%%%%%%%%%%%%%%%%%% END COROLLARY %%%%%
%
%%%%%%%%%%%%%%%%%%%%%%%%%%%%%%%%%%%%%%%%%%%%%%%%%%%%% PROOF %%%%%
\begin{proof}
Suppose that there exists a connected almost alternating diagram 
$\t{L}$ of the trivial link $L$ with $n$ ($> 2$) components. 
Here we assume that $\t{L}$ has the minimal crossings among such 
diagrams, and then  $\t{L}$ is prime. Since $L$ is also splittable, 
$\t{L}$ is not reduced from Theorem \ref{ThmLinks}. If $\t{L}$ is 
diagram II, we obtain another almost alternating link diagram of 
$L$ with less crossings than those of $\t{L}$ by applying 
reducing move II. This contradicts the minimality of the number 
of crossings of $\t{L}$. Therefore, $\t{L}$ is diagram I. Then 
apply reducing move I to $\t{L}$. If we obtain a connected 
diagram, then it is a connected alternating link diagram, which 
is non-splittable. Thus, we have a disconnected alternating 
diagram consisting of two connected components. Since $L$ has 
more than two components, at least one of them has more than 
one component, which is a connected alternating link diagram 
and thus non-splittable. Thus, $L$ is non-trivial, which 
contradicts our assumption.
\end{proof}\par
%%%%%%%%%%%%%%%%%%%%%%%%%%%%%%%%%%%%%%%%%%%%%%%%% END PROOF %%%%%
%
%%%%%%%%%%%%%%%%
\par\vspace{5mm}
%%%%%%%%%%%%%%%%
%
%
%
%%%%%%%%%%%%%%%%%%%%%%%%%%%%%%%%%%%%%%%%%%%%%%%%%%
\section{Standard Position for A Splitting Sphere} 
\label{SecStdPosition}
%%%%%%%%%%%%%%%%%%%%%%%%%%%%%%%%%%%%%%%%%%%%%%%%%%
%
%%%%%%%%%%%%%%%%%%%%%%%%%%%%%%%%%%%
Let $\t{L}$ $\subset$ $\bb{R}^2$ $\subset$ $S^2$ $=$ $\bb{R}^2$ 
$\cup$ $\{\infty\}$ be a diagram of a splittable link $L$. Here, 
we do not assume that $\t{L}$ is almost alternating. Note that 
we may speak sensibly about points ``above" or ``below" $\t{L}$ 
and also about ``inside" or ``outside" of some reasion, since we 
consider the projection plane $\bb{R}^2$ as a subspace of a 
$2$-sphere in $S^3$. 
%%%%%%%%%%%%%%%%%%%%%%%%%%%%%%%%%%%
%
%%%%%%%%%%%%%%%%
\par\vspace{3mm}
%%%%%%%%%%%%%%%%
%
%%%%%%%%%%%%%%%%%%%%%%%%%%%%%%%%%%%
Following \cite{m} and \cite{mt}, put a $3$-ball, called
a {\it crossing ball}, at each crossing point of $\t{L}$. Then,
isotope $L$ so that, at each crossing point, the overstrand runs
on the upper hemisphere and the understrand runs on the lower 
hemisphere as shown in Figure \ref{fig:crossingball}. We call 
the boundary of such a $3$-ball a {\it bubble}.
%%%%%%%%%%%%%%%%%%%%%%%%%%%%%%%%%%%
%
%%%%%%%%%%%%%%%%%%%%%%%%%%%%%%%%%%%
\par\vspace{3mm}\begin{figure}[htbp]\caption{A crossing ball}
\label{fig:crossingball}\end{figure}\par\vspace{3mm}
%%%%%%%%%%%%%%%%%%%%%%%%%%%%%%%%%%%
%
%%%%%%%%%%%%%%%%%%%%%%%%%%%%%%%%%%%
Let $S_+$ (resp. $S_-$) be the $2$-sphere obtained from $S^2$
by replacing the intersection disk of each crossing ball and 
$S^2$ by the upper (resp. lower) hemisphere of the bubble of 
the crossing ball. We will use the notation $S_{\pm}$ to mean 
$S_+$ or $S_-$ and similarly for other symbols with subscript 
$\pm$.
%%%%%%%%%%%%%%%%%%%%%%%%%%%%%%%%%%%
%
%%%%%%%%%%%%%%%%
\par\vspace{3mm}
%%%%%%%%%%%%%%%%
%
%%%%%%%%%%%%%%%%%%%%%%%%%%%%%%%%%%%
Let $F \in S^3 - L$ be a splitting sphere for a splittable link 
$L$. We may isotope $F$ to a suitable position with respect to 
$\t{L}$ according to \cite{m} and \cite{mt}.
%%%%%%%%%%%%%%%%%%%%%%%%%%%%%%%%%%%
%
%%%%%%%%%%%%%%%%%%%%%%%%%%%%%%%%%%% PROPOSITION PropIsotopy %%%%%
\begin{Prop}\label{PropIsotopy} Let $F$ be a surface mentioned 
above. Then, we may isotope $F$ so that;
%
%%%%%%%%%%%%%%%%%%%%%%%%%%%%%%%%%%%
\begin{itemize}
\item[(i)]  $F$ meets $S_{\pm}$ transversely in a pairwise 
            disjoint collection of simple closed curves and
\item[(ii)] $F$ meets each crossing-ball in a collection of 
            saddle-shaped disks $($Figure $\ref{fig:saddle}$$)$.
\end{itemize}
%%%%%%%%%%%%%%%%%%%%%%%%%%%%%%%%%%%
%
\end{Prop}
%%%%%%%%%%%%%%%%%%%%%%%%%%%%%%%%%%%%%%%%%%% END PROPOSITION %%%%%
%
%%%%%%%%%%%%%%%%%%%%%%%%%%%%%%%%%%%
\begin{figure}[htbp]\caption{A saddle}
\label{fig:saddle}\end{figure}\par\vspace{3mm}
%%%%%%%%%%%%%%%%%%%%%%%%%%%%%%%%%%%
%
%%%%%%%%%%%%%%%%%%%%%%%%%%%%%%%%%%%
Let $F$ be a splitting sphere satisfying the conditions in 
Proposition \ref{PropIsotopy}. The {\it complexity} $c(F)$ 
of $F$ is the lexicographically ordered pair $(t,u)$, where 
$t$ is the number of saddle-intersections of $F$ with 
crossing-balls of the diagram $\t{L}$, and $u$ is the total 
number of components of $F$ $\cap$ $S_+$ and $F$ $\cap$ $S_-$. 
We say that $F$ has {\it minimal complexity} if $c(F)$ $\leq$ 
$c(F')$ for any splitting sphere $F'$. Then, we have the 
following also according to \cite{m} and \cite{mt}.
%%%%%%%%%%%%%%%%%%%%%%%%%%%%%%%%%%%
%
%%%%%%%%%%%%%%%%%%%%%%%%%%%%%%%%%%% PROPOSITION PropMinimal %%%%%
\begin{Prop}\label{PropMinimal}
Let $F$ be a splitting sphere for a splittable link. If $F$ 
satisfies the conditions in Proposition \ref{PropIsotopy} 
and has minimal complexity, then each simple closed curve 
$C$ in $F$ $\cap$ $S_+$ $($resp. $F$ $\cap$ $S_-$$)$ meets 
the following requirements;
%%%%%%%%%%%%%%%%%%%%%%%%%%%%%%%%%%%
%
%%%%%%%%%%%%%%%%%%%%%%%%%%%%%%%%%%%
\begin{itemize}
\item[(i)]   $C$ bounds a disk in $F$ whose interior lies 
             entirely above $S_+$ $($resp. below $S_-$$)$,
\item[(ii)]  $C$ meets at least one bubble, and
\item[(iii)] $C$ does not meet any bubble in more than one arc.
\end{itemize}
%%%%%%%%%%%%%%%%%%%%%%%%%%%%%%%%%%%
%
\end{Prop}
%%%%%%%%%%%%%%%%%%%%%%%%%%%%%%%%%%%%%%%%%%% END PROPOSITION %%%%%
%
%%%%%%%%%%%%%%%%%%%%%%%%%%%%%%%%%%%
We say that a splitting sphere is in {\it standard position} 
if it satifies conclusions of Proposition \ref{PropIsotopy} 
and \ref{PropMinimal}. Now let $\t{L}$ be an almost alternating 
diagram of a splittable link $L$. Assume that $F$ is a splitting 
sphere for $L$ in standard position. Let $A$ denote the crossing 
ball at the dealternator. Note that our almost alternating 
diagram has exactly one dealternator. Let $C$ be a simple closed 
curve of $F \cap S_{\pm}$. Since $\t{L}$ is almost alternating, 
a subarc of $C$ satisfies the following {\it almost alternating 
property};
%%%%%%%%%%%%%%%%%%%%%%%%%%%%%%%%%%%
%
%%%%%%%%%%%%%%%%
\par\vspace{3mm}
%%%%%%%%%%%%%%%%
%
%%%%%%%%%%%%%%%%%%%%%%%%%%%%%%%%%%%
\noindent{\it If $C$ meets two bubbles of crossing-balls $B_1$ 
and $B_2$ in succession. Then,
%%%%%%%%%%%%%%%%%%%%%%%%%%%%%%%%%%%
%
%%%%%%%%%%%%%%%%%%%%%%%%%%%%%%%%%%%
\begin{itemize}
\item[(i)]  two arcs of $\t{L}$ $\cap$ $S_{\pm}$ on 
            $B_1$ and $B_2$ lie on the opposite sides of $C$ if 
            none of $B_1$ and $B_2$ are $A$ and
\item[(ii)] two arcs of $\t{L}$ $\cap$ $S_{\pm}$ on 
            $B_1$ and $B_2$ lie on the same side of $C$ if 
            one of $B_1$ and $B_2$ is $A$.
\end{itemize}}
%%%%%%%%%%%%%%%%%%%%%%%%%%%%%%%%%%%
%
%%%%%%%%%%%%%%%%
\par\vspace{3mm}
%%%%%%%%%%%%%%%%
%
%%%%%%%%%%%%%%%%%%%%%%%%%%%%%%%%%%%
Moreover, $C$ satisfies the following.
%%%%%%%%%%%%%%%%%%%%%%%%%%%%%%%%%%%
%
%%%%%%%%%%%%%%%%%%%%%%%%%%%%%%%%%%%%%%% LEMMA LemEveryCurve %%%%%
\begin{Lem}\label{LemEveryCurve} Every curve must pass the 
dealternator exactly once.   \end{Lem}
%%%%%%%%%%%%%%%%%%%%%%%%%%%%%%%%%%%%%%%%%%%%%%%%% END LEMMA %%%%%
%
%%%%%%%%%%%%%%%%%%%%%%%%%%%%%%%%%%%%%%%%%%%%%%%%%%%%% PROOF %%%%%
\begin{proof}
We may assume that the complexity of $F$ is finite.
It is sufficient to show just that $C$ passes the dealternator, 
since we have the third condition of Proposition 
\ref{PropMinimal}. Suppose that $C$ does not pass the 
dealternator. Here, we may assume that the dealternator is 
outside of $C$, that is, in the region with $\{\infty\}$ of the
two regions devided by $C$. Note that the length of every curve 
is more than one and even. From the almost alternating property, 
$C$ must contain at least one curve inside of it, say $C'$.
Then, $C'$ is also does not pass the dealternator, because 
it is inside of $C$. Since $C'$ also must contain at least 
one curve inside of it, we can inductively find 
an infinitely many curves inside of $C$. This contradicts 
the finiteness of the complexity of $F$.
\end{proof}\par
%%%%%%%%%%%%%%%%%%%%%%%%%%%%%%%%%%%%%%%%%%%%%%%%% END PROOF %%%%%
%
%%%%%%%%%%%%%%%%%%%%%%%%%%%%%%%%%%%
The collection of circles of $F$ $\cap$ $S_{\pm}$, together with 
the saddle components of $F$ $\cap$ $($$\cup$ $\{B_i\}$$)$, give 
rise to a cell-decomposition of $F$, which we call the
intersection graph $G$ of $F$ with $S_{\pm}$ and $\cup$ 
$\{B_i\}$. The vertices of $G$ correspond to the saddles of $F$ 
$\cap$ $($$\cup$ $\{B_i\}$$)$, the edges of $G$ correspond to 
the arcs of $F$ $\cap$ $($$S_+$ $-$ $\cup$ $\{B_i\}$$)$ $=$ 
$F$ $\cap$ $($$S_-$ $\cup$ $\{B_i\}$$)$, and the faces of $G$ 
correspond to the disks, called {\it dome}, bounded by the 
simple closed curves of $F$ $\cap$ $S_+$ and of $F$ $\cap$ $S_-$, 
afforded to us by the first condition of Proposition 
\ref{PropMinimal}. Note that $G$ is a plane graph in sphere $F$ 
and the degree of each vertex of $G$ is $4$. We define the degree 
of a face as the number of vertices which the face has on its 
boundary. Let $f_i$ and $\mid f_i \mid$ be a face with degree $i$ 
and the number of faces of degree $i$, respectively. Then, we have 
the following from the Euler's formula.
%%%%%%%%%%%%%%%%%%%%%%%%%%%%%%%%%%%
%
%%%%%%%%%%%%%%%%%%%%%%%%%%%%%%%%%%%%%%%%%%%% LEMMA LemEuler %%%%%
\begin{Lem}\label{LemEuler}  $\sum (i-4) \mid f_i \mid = -8$
\end{Lem}
%%%%%%%%%%%%%%%%%%%%%%%%%%%%%%%%%%%%%%%%%%%%%%%%% END LEMMA %%%%%
%
%%%%%%%%%%%%%%%%%%%%%%%%%%%%%%%%%%%%%%%%%%%%%%%%%%%%% PROOF %%%%%
\begin{proof}
Let $n$, $e$, and $f$ be the numbers of the vertices, edges, and 
faces of $G$, respectively. From the Euler's formula, we have 
$n$ $-$ $e$ $+$ $f$ $=$ $2$. Since we also have $\sum$ $i$ 
$\mid f_i \mid$ $=$ $4n$, $\sum$ $i$ $\mid f_i \mid$ $=$ $2e$ 
and $\sum$ $\mid f_i \mid$ $=$ $f$, $\sum$ $(i-4)$ 
$\mid f_i \mid$ $=$ $-$ $\sum$ $i$ $\mid f_i \mid$ $+$ $2$ 
$\sum$ $i$ $\mid f_i \mid$ $-$ $4$ $\sum$ $\mid f_i \mid$ 
$=$ $-$ $4n$ $+$ $4e$ $-$ $4f$ $=$ $-8$.
\end{proof}\par
%%%%%%%%%%%%%%%%%%%%%%%%%%%%%%%%%%%%%%%%%%%%%%%%% END PROOF %%%%%
%
%%%%%%%%%%%%%%%%%%%%%%%%%%%%%%%%%%%
For a convenience, we introduce several terminologies. We call a 
vertex a {\it black vertex} and denote it by $v_d$ if it comes 
from a saddle on the dealternator. And also, we denote by $b_d$ 
the bubble put on the dealternator. Moreover, denote by $b_i$ the 
bubble which contains a saddle corresponding to vertex $v_i$, 
which we call a {\it white vertex}. Put color black and white to 
black vertices and to other vertices, respectively. We say, as 
usual, a face is {\it adjacent} to another face if they have a 
common edge on their boundaries. We say that a face is {\it 
vertexwise-adjacent} to another face if they have a common vertex 
on their boundaries.                
%%%%%%%%%%%%%%%%%%%%%%%%%%%%%%%%%%%
%
%%%%%%%%%%%%%%%%
\par\vspace{5mm}
%%%%%%%%%%%%%%%%
%                                                        
%                                                                
%
%                                              
%%%%%%%%%%%%%%%%%%%%%%%%%%%%%%%%%%%%%%%%%                   
\section{Proof of Theorem \ref{ThmLinks}} \label{SecPrfLinks}   
%%%%%%%%%%%%%%%%%%%%%%%%%%%%%%%%%%%%%%%%%
%
%%%%%%%%%%%%%%%%%%%%%%%%%%%%%%%%%%%  
We need the following lemma to prove Theorem \ref{ThmLinks}.     
Remark here that our almost alternating diagram has exactly 
one dealternator.     
%%%%%%%%%%%%%%%%%%%%%%%%%%%%%%%%%%%  
%                                                                
%
%
%%%%%%%%%%%%%%%%%%%%%%%%%%%%%%%%%%%%%%% LEMMA LemMinimality %%%%%
\begin{Lem}\label{LemMinimality} Let $\t{L}$ be a connected, 
prime, reduced almost alternating diagram of a splittable 
link $L$. If $\t{L}$ is one of the diagrams in Figure 
\ref{fig:NotMinimal}, then there exists another connected, 
prime, reduced almost alternating diagram of $L$ or of 
another splittable link with less crossings than $\t{L}$.    
\end{Lem}                                                       
%%%%%%%%%%%%%%%%%%%%%%%%%%%%%%%%%%%%%%%%%%%%%%%%% END LEMMA %%%%%
%         
%%%%%%%%%%%%%%%%%%%%%%%%%%%%%%%%%%%  
\begin{figure}[htbp]\caption{}
\label{fig:NotMinimal}\end{figure}\par\vspace{3mm}
%%%%%%%%%%%%%%%%%%%%%%%%%%%%%%%%%%%  
%    
%%%%%%%%%%%%%%%%%%%%%%%%%%%%%%%%%%%%%%%%%%%%%%%%%%%%% PROOF %%%%%
\begin{proof}
We show only the case for diagram VI, which has three tangle 
areas $T_1$, $T_2$, and $T_3$ and ten regions $a$, $b$, $\ldots$, 
$j$. Some of these regions might be the same. Then we have 
six possibilities which regions are the same from the reducedness 
and the primeness of $\t{L}$ (for instance, we have a nonprime 
diagram if $a$ $=$ $d$ and we have a nonreduced diagram if $b$ 
$=$ $e$ and $c$ $=$ $f$). Here, we show the case that all ten 
regions are mutually different. Other cases can be shown 
similarly. In addition, let us assume that regions $g$ and $j$
do not share an arc inside of $T_3$ (in this case, we consider
diagram $\t{K}$ of $L$ and diagram $\t{K'}$ of another splittable
link instead of $\t{L}$ and $\t{L'}$, see Figure \ref{fig:KandKP}).
If link $L$ has diagram $\t{L}$, then $L$ has diagram $\t{L'}$ 
with one less crossings than those of $\t{L}$ as well. Then, note 
that all nine regions $k$, $l$, $\ldots$, $s$ are mutually different.
%%%%%%%%%%%%%%%%%%%%%%%%%%%%%%%%%%%
%
%%%%%%%%%%%%%%%%%%%%%%%%%%%%%%%%%%%  
\begin{figure}[htbp]\caption{$\t{L}$ and $\t{L'}$}
\label{fig:LandLP}\end{figure}\par\vspace{3mm}
\begin{figure}[htbp]\caption{$\t{K}$ and $\t{K'}$}
\label{fig:KandKP}\end{figure}\par\vspace{3mm}
%%%%%%%%%%%%%%%%%%%%%%%%%%%%%%%%%%%  
%
%%%%%%%%%%%%%%%%%%%%%%%%%%%%%%%%%%% 
{\bf (Connectedness)} \hspace{3mm} Assume that $\t{L'}$ is not
connected. Then, we have a simple closed surve $C$ in a region
of $\t{L'}$ such that each of the two regions of $S^2$ $-$ $C$ 
contains a component of $\t{L'}$ (here we call such a curve a
{\it splitting curve}). If $C$ is entirely contained in a 
tangle area, then it is easy to see that $\t{L}$ is not 
connected as well. Therefore, $C$ $\cap$ $\{$ $\cup$ $T_i$ $\}$ 
$\ne$ $\emptyset$ and $C$ is in region $k$, $l$, $\ldots$, or 
$s$. We may assume that $C$ has minimal intersection with 
$\cup$ $T_i$. Take a look at one of outermost intersections of
$T_i$ and $C$. Then the intersection is one of the following 
(Figure \ref{fig:ConIntscn}). In the cases of (i) and (iv), we 
can have another splitting curve $C'$ which is entirely 
contained in $T_i$ or which has one less intersections with 
$\cup$ $T_i$ than $C$ does. In the cases of (iii) and (vi), 
we have the same regions among the nine regions of $\t{L'}$. In 
other cases, $C$ has an intersection with $\t{L'}$.
%%%%%%%%%%%%%%%%%%%%%%%%%%%%%%%%%%%
%         
%%%%%%%%%%%%%%%%%%%%%%%%%%%%%%%%%%%  
\begin{figure}[htbp]\caption{}
\label{fig:ConIntscn}\end{figure}\par\vspace{3mm}
%%%%%%%%%%%%%%%%%%%%%%%%%%%%%%%%%%%  
%    
%%%%%%%%%%%%%%%%
\par\vspace{3mm}
%%%%%%%%%%%%%%%%
%
%%%%%%%%%%%%%%%%%%%%%%%%%%%%%%%%%%% 
{\bf (Primeness)} \hspace{3mm} Assume that $\t{L'}$ is not
prime. Then we have a simple closed curve $C$ which intersects
$\t{L'}$ in just two points belonging to different arcs of 
$\t{L'}$ (here we call such a curve a {\it separating curve}).
If $C$ is in a tangle area, then it is easy to see that $\t{L}$ 
is not prime as well. However, $C$ $\cap$ $\{$ $\cup$ $T_i$ 
$\}$ $\ne$ $\emptyset$, since no pair of nine regions share two 
different arcs outside of tangle areas. We may assume that 
$C$ has minimal intersection with $\cup$ $T_i$. Take a look at 
one of the outermost intersections of $T_i$ and $C$. Then the 
intersection is one of the figures in Figure \ref{fig:ConIntscn}.
In the cases (i) and (iv), we can eliminate the intersection,
which is a contradiction. In the cases of (ii) and (v), we can
obtain another separating curve $C'$ which is entirely contained 
in $T_i$ or which has one less intersections with $\cup$ $T_i$ 
than $C$ does. In the cases of (iii) and (vi), we have the same 
regions among the nine regions of $\t{L'}$. In the case of 
(vii), $C$ has an intersection with $T_3$ and here we may assume 
that the other intersection is outside of $T_3$ from 
the minimality. Then regions $k$ and $n$ must share arcs inside
and outside of $T_3$, and thus regions $g$ and $j$ must share an
arc inside of $T_3$, which contradicts our assumption.
%%%%%%%%%%%%%%%%%%%%%%%%%%%%%%%%%%%
%            
%%%%%%%%%%%%%%%%
\par\vspace{3mm}
%%%%%%%%%%%%%%%%
%
%%%%%%%%%%%%%%%%%%%%%%%%%%%%%%%%%%% 
{\bf (Reducedness)} If $\t{L'}$ is diagram I, then we obtain an 
alternating diagram of $L$ by reducing move I. Thus, $L$ is 
non-splittable, which is a contradiction. Assume that $\t{L'}$ 
is diagram II. Then we can find a part in the diagram which we 
can apply reducing move II to (Figure \ref{fig:NonRedPart}). We 
have four possibilities; $(x, y)$ $=$ $(n,o)$, $(o,p)$, $(p,q)$, 
or $(q,n)$. In the first case, regions $q$ and $u$ must share a 
crossing so that $\t{L'}$ contains the part in Figure 
\ref{fig:NonRedPart}. Since $m$ $\ne$ $o$ and $l$ $\ne$ $n$, we 
have regions $t$ and $u$, and then $u$ might be the same as $m$ 
or $o$ (see Figure \ref{fig:NonRedCsOne}). However, the regions 
which can share a crossing with region $q$ are $k$, $o$, or $s$. 
If $u$ $=$ $o$, then we obtain a non-prime diagram. Therefore, 
this case does not occure.
%%%%%%%%%%%%%%%%%%%%%%%%%%%%%%%%%%%
%         
%%%%%%%%%%%%%%%%%%%%%%%%%%%%%%%%%%% 
In the second case, regions $k$ and $q$ must share a crossing.
Thus, we can decide the inside of $T_3$ more precisely and
then we can see that $\t{L}$ is non-reduced (Figure 
\ref{fig:NonRedCsTwo}).
%%%%%%%%%%%%%%%%%%%%%%%%%%%%%%%%%%%
%         
%%%%%%%%%%%%%%%%%%%%%%%%%%%%%%%%%%% 
In the third case, regions $o$ and $v$ must share a crossing
(Figure \ref{fig:NonRedCsThree}). The regions which can share
a crossing with $o$ are $k$, $m$, or $q$. We also have three
pssibilities that $v$ $=$ $k$, $q$, or $s$. Thus, we obtain
that $v$ $=$ $q$ or $k$. In the former case, we have a non-prime
diagram. In the latter case, we can decide the inside of $T_3$
more precisely and then we can see that $\t{L}$ is non-reduced
(Figure \ref{fig:NonRedCsThree}).
%%%%%%%%%%%%%%%%%%%%%%%%%%%%%%%%%%%
%         
%%%%%%%%%%%%%%%%%%%%%%%%%%%%%%%%%%%
We can prove the fourth case similarly.
\end{proof}   
%%%%%%%%%%%%%%%%%%%%%%%%%%%%%%%%%%%  
%      
%%%%%%%%%%%%%%%%%%%%%%%%%%%%%%%%%%%  
\begin{figure}[htbp]\caption{}
\label{fig:NonRedPart}\end{figure}\par\vspace{3mm}
\begin{figure}[htbp]\caption{}
\label{fig:NonRedCsOne}\end{figure}\par\vspace{3mm}
\begin{figure}[htbp]\caption{}
\label{fig:NonRedCsTwo}\end{figure}\par\vspace{3mm}
\begin{figure}[htbp]\caption{}
\label{fig:NonRedCsThree}\end{figure}\par\vspace{3mm}
%%%%%%%%%%%%%%%%%%%%%%%%%%%%%%%%%%% 
% 
%%%%%%%%%%%%%%%%%%%%%%%%%%%%%%%%%%%%%%%%%%%%%%%%% END PROOF %%%%% 
%
%%%%%%%%%%%%%%%%%%%%%%%%%%%%%%%%%%%%%%%%%%%%%%%%%%%%%%%%%%%%%%%%%
\noindent{\bf Proof of Theorem \ref{ThmLinks}}.                 
%%%%%%%%%%%%%%%%%%%%%%%%%%%%%%%%%%%%%%%%%%%%%%%%%%%%%%%%%%%%%%%%%
%                                                                
%%%%%%%%%%%%%%%%%%%%%%%%%%%%%%%%%%%  
Suppose that there exists a splittable link with a connected, 
prime, reduced almost alternating diagram. Take all such links 
and consider all such diagrams of them. Let $\t{L}$ be minimal 
in such diagrams with respect to the number of crossings. Then 
$\t{L}$ is none of the diagrams in Figure \ref{fig:NotMinimal}, 
otherwise it contradicts the minimality of $\t{L}$ from Lemma 
\ref{LemMinimality}. Note that $\t{L}$ has at least two 
crossings, since $\t{L}$ is connected and $\t{L}$ has more than 
one component. Let $F$ $\subset$ $S^3$ $-$ $L$ be a splitting 
sphere for $L$, which would be assumed to be in a standard 
position. And let $G$ $\subset$ $F$ be the intersection graph of 
$F$ $\cap$ $S^2$. For each face of degree $i$ of $G$, charge 
weight $i-4$. We denote by $w(f)$ the weight of a face $f$. If 
there is no faces of degree $2$, then every face has non negative 
weight. Then, $\sum$ $(i-4)$ $\mid f_i \mid$ $\geq$ $0$ (the sum 
of weights of all faces), which contradicts Lemma \ref{LemEuler}. 
Therefore, we may assume that there exists at least one face of 
degree $2$.        
%%%%%%%%%%%%%%%%%%%%%%%%%%%%%%%%%%%
%
%%%%%%%%%%%%%%%%
\par\vspace{3mm}
%%%%%%%%%%%%%%%%
%
%%%%%%%%%%%%%%%%%%%%%%%%%%%%%%%%%%%                           
It may happen that two faces of degree $2$ are adjacent or 
vertexwise-adjacent to each other. However if two faces of 
degree $2$ are adjacent to each other, then it contradicts 
the reducedness of $\t{L}$ (Figure \ref{fig:fook}). Also if 
two faces of degree $2$ are vertexwise-adjacent to each other 
at a white vertex, then there exists a face which has two 
black vertices on its boundary, which contradicts Lemma 
\ref{LemEveryCurve}. Therefore, we have two cases if we look 
at a face of degree $2$. One is that it is not adjacent or 
vertexwise-adjacent to any other faces of degree $2$. Here we 
call it a block of type $T'$ or simply $T'$ and then $w(T') =  
-2$. The other is that it is vertexwise-adjacent to another 
face of degree $2$ at a black vertex. In this case, we put 
these two faces together and call it a block of type $U'$ or 
simply $U'$, and then $w(U') = -4$, which is the sum of the 
weights of the two faces of degree $2$. 
%%%%%%%%%%%%%%%%%%%%%%%%%%%%%%%%%%%
%
%%%%%%%%%%%%%%%%%%%%%%%%%%%%%%%%%%%
\par\vspace{3mm}\begin{figure}[htbp]\caption{}
\label{fig:fook}\end{figure}\par\vspace{3mm}
%%%%%%%%%%%%%%%%%%%%%%%%%%%%%%%%%%%
%     
%%%%%%%%%%%%%%%%%%%%%%%%%%%%%%%%%%%
Take a look at two faces $f_{i\ge4}$ and $f_{j\ge4}$ which are
adjacent to a block of type $T'$ (resp. $U'$) and put all of them 
together. We call it a block of type $Y_{i,j}$ (resp. $Z_{i,j}$)
or simply $Y_{i,j}$ (resp. $Z_{i,j}$). In the case of $Z_{i,j}$,
we assume that $i$ is greater than equal to $j$. Then, 
$w(Y_{i,j})$ $=$ $w(f_i)$ $+$ $w(f_j)$ $+$ $w(T')$ $=$ 
$i$ $+$ $j$ $-$ $10$ $\ge$ $-2$ and 
$w(Z_{i,j})$ $=$ $w(f_i)$ $+$ $w(f_j)$ $+$ $w(U')$ $=$ 
$i$ $+$ $j$ $-$ $12$ $\ge$ $-4$. We call a face of degree $i$
($\ge4$) a block of type $X_i$ or simply $X_i$ if it is not 
adjacent to $T'$ or $U'$. Then, $w(X_i)$ $=$ $i$ $-$ $4$ $\ge$ 
$0$. For blocks of type $Y_{i,j}$, the only type of blocks with
a negative weight is $Y_{4,4}$ and we call a block of type
$Y_{4,4}$ a block of type $T$ or simply $T$. If there exists
$Z_{4,4}$, then it contradicts the minimality of $\t{L}$, since
we have diagram III (Figure \ref{fig:z44}). Thus, for 
blocks of type $Z_{i,j}$, the only type of blocks with a 
negative weight is $Z_{6,4}$ and we call a block of type 
$Z_{6,4}$ a block of type $U$ or simply $U$.
%%%%%%%%%%%%%%%%%%%%%%%%%%%%%%%%%%%
%
%%%%%%%%%%%%%%%%%%%%%%%%%%%%%%%%%%%
\par\vspace{3mm}\begin{figure}[htbp]\caption{$Z_{4,4}$}
\label{fig:z44}\end{figure}\par\vspace{3mm}
%%%%%%%%%%%%%%%%%%%%%%%%%%%%%%%%%%%
% 
%%%%%%%%%%%%%%%%%%%%%%%%%%%%%%%%%%%
If there are no blocks of type $T$ nor type $U$, then it 
contradicts Lemma \ref{LemEuler} as before. Here we say that a 
block is {\it upper} (resp. {\it lower}) if its faces ($\ne$ 
$f_2$) come from domes which are above $S_+$ (resp. below $S_-$). 
Consider the following three cases; $G$ has $T_+$ and no $U_+$ 
({\bf Case 1}), $G$ has $U_+$ and no $T_+$ ({\bf Case 2}), and 
$G$ has $T_+$ and $U_+$ ({\bf Case 3}), where $T_+$ means an 
upper block of type $T$ and $U_-$ means a lower block of type 
$U$, for instance. In each case, we show that we can discharge 
weights of lower blocks to $T_+$ and $U_+$ to make the weight 
of every block non-negative. Therefore, proving the above three 
cases tells us that there does not exist graph $G$, that is, 
there does not exist a connected, prime reduced almost 
alternating diagram of any splittable link. This completes the 
proof.
%%%%%%%%%%%%%%%%%%%%%%%%%%%%%%%%%%%
%
%%%%%%%%%%%%%%%%
\par\vspace{3mm}
%%%%%%%%%%%%%%%%
%
%%%%%%%%%%%%%%%%%%%%%%%%%%%%%%%%%%%
In each  of three cases, we induce a contradiction by actually 
replacing the boundary cycles of subgraphs of $G$ on the diagram 
and looking at the diagram as shown in Figure \ref{fig:fook} 
or Figure \ref{fig:z44}. Here, put orientations on $S^2$ and $F$.
We have two possibilities to replace the boundary cycle of a face
on the diagram; its orientation coincides that of $S^2$ or not
(we did not mention about this before). However, we may 
occasionaly choose one of the two possibilities, since the 
diagrams obtained by the two ways are the same up to mirror image, 
which does not affect our purpose. Here we have the following 
claim.
%%%%%%%%%%%%%%%%%%%%%%%%%%%%%%%%%%%
%
%%%%%%%%%%%%%%%%%%%%%%%%%%%%%%%%%%%%%%%%%%%%% CLAIM ClmBdTU %%%%%
\begin{Claim}\label{ClmBdTU}\begin{itemize}
\item[(i)]  For every block of type $T$, its five white vertices 
            come from saddles in mutually different five bubbles.
\item[(ii)] For blocks of type $U$, the boundary curves of the
            faces of degree four pass the same four bubbles.
\end{itemize}\end{Claim}\par    
%%%%%%%%%%%%%%%%%%%%%%%%%%%%%%%%%%%%%%%%%%%%%%%%% END CLAIM %%%%%
%
%%%%%%%%%%%%%%%%%%%%%%%%%%%%%%%%%%%%%%%%%%%%%%%%%%%%% PROOF %%%%%
\begin{proof}                                                   
(i) Take a block of type $T$ and put names $v_\a$, $v_\b$, $v_\c$, 
$v_\d$, and $v_\e$ to it as shown in Figure \ref{fig:TandU}. 
Assume that there is a pair of vertices coming from saddles in 
a same bubble. From the almost alternating property, we have that 
$b_\a$ $=$ $b_\e$, $b_\a$ $=$ $b_\d$, or $b_\b$ $=$ $b_\e$. The 
first case contradicts the mimimality of $\t{L}$ (diagram III)
and the last two cases contradict the primeness.  
(ii) Take two blocks of type $U$. Put names $v_\f$, $v_\g$, and 
$v_\h$ to one of them and $v_{\f'}$, $v_{\g'}$, and $v_{\h'}$ to
the other following Figure \ref{fig:TandU}. Then, we have that 
$b_{\f'}$ $=$ $b_\h$ and $b_{\h'}$ $=$ $b_\f$ or that $b_{\f'}$ 
$=$ $b_\f$ and $b_{\h'}$ $=$ $b_{\h}$. The former case 
contradicts the minimality of $\t{L}$ (diagram III) and the 
latter case contradicts the primeness unless the claim holds.     
\end{proof} 
%%%%%%%%%%%%%%%%%%%%%%%%%%%%%%%%%%%%%%%%%%%%%%%%% END PROOF %%%%%
%
%%%%%%%%%%%%%%%%%%%%%%%%%%%%%%%%%%%
\par\vspace{3mm}\begin{figure}[htbp]\caption{T and U}
\label{fig:TandU}\end{figure}\par\vspace{3mm}
%%%%%%%%%%%%%%%%%%%%%%%%%%%%%%%%%%%
% 
%%%%%%%%%%%%%%%%%%%%%%%%%%%%%%%%%%%
Since the boundary curves of the faces of degree $2$ (resp. 
$4$) of all blocks of type $T_{\pm}$ (resp. $U_{\pm}$) pass 
the same two (resp. four) bubbles and are parallel (otherwise,
it contradicts the primeness or the reducedness), we can define 
above, below, the leftside of, and the rightside of the 
dealternator on the diagram as shown in Figure \ref{fig:ABLR}. 
We define the top and the bottom face of $T$ (resp. the left 
and the right face of $U$) as the face of degree $4$ (resp. $2$) 
which is above and below (resp. the leftside of and the 
rightside of) the dealternator on the diagram, respectively. 
To the boundary curves of two faces which are not 
vertexwise-adjacent to each other at the dealternator, we 
define that one is outside of the other if it is closer to the 
center of the dealternator than the other is on the diagram (see 
Figure \ref{fig:IO}). Before we start, we define the following 
three types of adjacency.
%%%%%%%%%%%%%%%%%%%%%%%%%%%%%%%%%%%
%
%%%%%%%%%%%%%%%%
\par\vspace{3mm}
%%%%%%%%%%%%%%%%
%
%%%%%%%%%%%%%%%%%%%%%%%%%%%%%%%%%%%%%%%%%%%%%%%%%%%%%%%%%%%%%%%%%
\begin{itemize}                                                  
 \item[(A)] If a face is vertexwise-adjacent to the face of 
            degree $2$ of a block of type $T$ at the white 
            vertex, then we say that the face is $A$-adjacent 
            to the block of type $T$.  
 \item[(B)] If a face is adjacent to the top (resp. the bottom) 
            face of a block of type $T$ at edge $v_\d$$v_\e$ 
            (resp. $v_\a$$v_\b$), then we say that the face is 
            $B_t$- (resp. $B_b$-) adjacent to the block of type 
            $T$.
 \item[(C)] If a face is vertexwise-adjacent to the left (resp.
            the right) face of a block of type $U$, then we say 
            that the face is $C_l$- (resp. $C_r$-) adjacent to 
            the block of type $U$. 
\end{itemize}                                                    
%%%%%%%%%%%%%%%%%%%%%%%%%%%%%%%%%%%%%%%%%%%%%%%%%%%%%%%%%%%%%%%%%
%
%%%%%%%%%%%%%%%%%%%%%%%%%%%%%%%%%%%
\par\vspace{3mm}\begin{figure}[htbp]\caption{}
\label{fig:ABLR}\end{figure}\par\vspace{3mm}
%%%%%%%%%%%%%%%%%%%%%%%%%%%%%%%%%%%
%
%%%%%%%%%%%%%%%%%%%%%%%%%%%%%%%%%%%
\begin{figure}[htbp]\caption{}
\label{fig:IO}\end{figure}\par\vspace{5mm}
%%%%%%%%%%%%%%%%%%%%%%%%%%%%%%%%%%%
%    
%
%%%%%%%%%%%%%%%%%%%%%
%
%%%%%%%%%%%%%%%%%%%%%%%%%%%%%%%%%%%%%%%%%%%%%%%%%%%% Case 1 %%%%%
%
%%%%%%%%%%%%%%%%%
%
%%%%%%%%%%%%%%%%%%%%%%%%%%%%%%%%%%%     
\begin{center}{\bf Case 1.}\end{center}
%%%%%%%%%%%%%%%%%%%%%%%%%%%%%%%%%%%     
%
%%%%%%%%%%%%%%%%
\par\vspace{5mm}
%%%%%%%%%%%%%%%%
%             
%%%%%%%%%%%%%%%%%%%%%%%%%%%%%%%%%%%     
We first look at faces which are $A$-adjacent to blocks of type 
$T$. Since faces of degree $2$ of all blocks of type $T$ pass 
the same two bubbles, every face can be $A$-adjacent to at most 
one $T$ at most once. Then we have $7$ types of blocks which are
$A$-adjacent to blocks of type $T$; $X_i^a$, $Y_{i,j}^{p,q}$, and
$Z_{i,j}^{p,q}$ with $\{p,q\}$ $=$ $\{\cdot,a\}$, $\{a,\cdot\}$,
or $\{a,a\}$, where $Z_{i,j}^{a,\cdot}$ stands for a block of
type $Z_{i.j}$ whose $f_i$ is $A$-adjacent to a block of type $T$
and $f_j$ is not, for instance. We generally use $v_\a$, $v_\b$, 
$v_\c$, $v_\d$, and $v_\e$ to represent vertices of a block of 
type $T$ as the proof of Claim \ref{ClmBdTU} (i), which ensures 
us that $b_\a$, $b_\b$, $b_\c$, $b_\d$, and $b_\e$ are mutually 
different.
%%%%%%%%%%%%%%%%%%%%%%%%%%%%%%%%%%%
%
%%%%%%%%%%%%%%%%
\par\vspace{3mm}
%%%%%%%%%%%%%%%%
%
%%%%%%%%%%%%%%%%%%%%%%%%%%%%%%%%%%%%%%%%%%%%% CLAIM ClmAadj %%%%%
\begin{Claim}\label{ClmAadj}\begin{itemize}                     
\item[(i)]   No face of degree $4$ can be adjacent to a face 
             of degree $2$ and $A$-adjacent to a block of 
             type $T$.  
\item[(ii)]  No face of degree $4$ can be adjacent to two 
             faces of degree $2$ with any face which is 
             $A$-adjacent to a block of type $T$.      
\item[(iii)] No face of degree $6$ can be adjacent to two 
             faces of degree $2$ with any other face of degree 
             $6$ which is $A$-adjacent to a block of type $T$. 
\end{itemize}\end{Claim}                                        
%%%%%%%%%%%%%%%%%%%%%%%%%%%%%%%%%%%%%%%%%%%%%%%%% END CLAIM %%%%%
%
%%%%%%%%%%%%%%%%%%%%%%%%%%%%%%%%%%%%%%%%%%%%%%%%%%%%% PROOF %%%%%
\begin{proof}                                                   
(i) It contradicts the mimimality of $\t{L}$ (diagram IV).
(ii) Assume that the boundary curves of $T$ and the face $f$ 
which is $A$-adjacent to it have been replaced on the diagram. 
Put names $v_1$, $v_2$, and $v_3$ to the face of degree $4$ as 
shown in Figure \ref{fig:z4z66}. We may assume that $b_1$ is 
below the dealternator. Here note that $v_1$ is on the boundary 
cycle of $f$. Therefore we have that $b_1$ is surrounded by 
boundary curve $b_\a$$b_\b$$b_\c$$b_d$ or that $v_1$ $=$ $v_\b$ 
and $v_2$ $=$ $v_\a$. Similarly, we have that $b_3$ is 
surrounded by boundary curve $b_\c$$b_\d$$b_\e$$b_d$ or that
$v_3$ $=$ $v_\d$ and $v_2$ $=$ $v_\e$. It is easy to see that
none of four cases can be held considering the length of the 
boundary curve of $f_4$. 
(iii) Put names $v_1$, $\ldots$, $v_8$ as shown in Figure 
\ref{fig:z4z66}. From symmetricity, we may assume that the 
face, say $f$, with $v_d$$v_1$$v_2$$v_3$$v_4$$v_5$ as its 
boundary cycle is $A$-adjacent to $T$ at $v_2$ or $v_3$. 
Assume that the boundary curves of $f$ and $T$ have been
replaced on the diagram. In the first case, we have that
$b_6$ $=$ $b_\a$, $b_6$ $=$ $b_\b$, $b_7$ $=$ $b_\b$, or 
$b_7$ $=$ $b_\a$ considering replacing the boundary cycle 
$v_d$$v_5$$v_6$$v_7$$v_8$$v_1$ on the diagram. The first 
three cases contradicts the mimimality of $\t{L}$ (diagram IV or
V) and the last case contradicts the primeness. 
In the second case, we have that $b_6$, $b_7$, or $b_8$ $=$ 
$b_\c$ considering replacing boundary cycle 
$v_d$$v_5$$v_6$$v_7$$v_8$$v_1$ on the diagram. The first and 
the third cases contradict the mimimality of $\t{L}$ (diagram IV)
and the second case contradicts the primeness.
\end{proof}   
%%%%%%%%%%%%%%%%%%%%%%%%%%%%%%%%%%%%%%%%%%%%%%%%% END PROOF %%%%%
%
%%%%%%%%%%%%%%%%%%%%%%%%%%%%%%%%%%%
\par\vspace{3mm}\begin{figure}[htbp]\caption{$Z_4$ and $Z_{6,6}$}
\label{fig:z4z66}\end{figure}\par\vspace{3mm}
%%%%%%%%%%%%%%%%%%%%%%%%%%%%%%%%%%%
%
%%%%%%%%%%%%%%%%%%%%%%%%%%%%%%%%%%%
We have that $w(X_i^a)$ $=$ $i-4$ $\ge$ $4-4$ $=$ $0$. From Claim
\ref{ClmAadj}, we obtain that $Y_{i,j}^{a,\cdot}$ $=$ 
$Y_{\ge6,\ge4}$,
$Y_{i,j}^{\cdot,a}$ $=$ $Y_{\ge4,\ge6}$,
$Y_{i,j}^{a,a}$ $=$ $Y_{\ge6,\ge6}$,
$Z_{i,j}^{a,\cdot}$ $=$ $Z_{\ge8,\ge6}$,
$Z_{i,j}^{\cdot,a}$ $=$ $Z_{\ge8,\ge6}$, and
$Z_{i,j}^{a,a}$ $=$ $Z_{\ge8,\ge6}$. Therefore, we have that
$w(Y_{i,j}^{a,\cdot})$ $=$ $i$ $+$ $j$ $-$ $10$ $\ge$ $0$,
$w(Y_{i,j}^{\cdot,a})$ $\ge$ $0$, 
$w(Y_{i,j}^{a,a})$ $\ge$ $2$,
$w(Z_{i,j}^{a,\cdot})$ $=$ $i$ $+$ $j$ $-$ $12$ $\ge$ $2$,
$w(Z_{i,j}^{\cdot,a})$ $\ge$ $2$, and
$w(Z_{i,j}^{a,a})$ $\ge$ $2$.
%%%%%%%%%%%%%%%%%%%%%%%%%%%%%%%%%%%
%
%%%%%%%%%%%%%%%%%%%%%%%%%%%%%%%%%%%
Then, for each block which is $A$-adjacent to blocks of type $T$, 
discharge $2$ out of its weight to each of the blocks of type $T$ 
if the sum of the weights of the block and all the blocks of type 
$T$ is non-negative. If the sum is negative, call it a block of 
type $\A$, $\B^*$, $\C^*$, $\D^*$, $\E^*$, $\F^*$, $\G^*$, $\H$, 
$\I^*$ or $\J^*$ as follows, where $\B^*$ means $\B$ or $\B'$, 
for instance.      
%%%%%%%%%%%%%%%%%%%%%%%%%%%%%%%%%%%
%
%%%%%%%%%%%%%%%%
\par\vspace{3mm}
%%%%%%%%%%%%%%%%
%
%%%%%%%%%%%%%%%%%%%%%%%%%%%%%%%%%%%
The type of a block such that the sum of the weights of the 
block and blocks of type $T$ which the block is $A$-adjacent to 
is negative is $X_4^a$, $Y_{6,4}^{a,\cdot}$, $Y_{4,6}^{\cdot,a}$, 
$Y_{6,6}^{a,a}$, or $Z_{8,6}^{a,a}$. We consider the first and 
the last three cases. In the second case, we obtain the same 
types as those of the third case. It is easy to see that we can 
uniquely obtain the diagram from $X_4^a$ with $T$ on the diagram 
and we say that the block has type $\A$.
%%%%%%%%%%%%%%%%%%%%%%%%%%%%%%%%%%%
%
%%%%%%%%%%%%%%%%
\par\vspace{3mm}
%%%%%%%%%%%%%%%%
%
%%%%%%%%%%%%%%%%%%%%%%%%%%%%%%%%%%%
Take a look at $Y_{4,6}^{\cdot,a}$ and put names $v_1$, 
$\ldots$, $v_7$ as shown in Figure \ref{fig:y46y66z86}. 
Then, its $f_6$ is $A$-adjacent to $T$ at $v_2$, $v_3$, 
or $v_4$. In each case, replace the boundary cycle of 
its $f_4$ on the diagram assuming that we have already 
replaced the boundary cycles of $f_6$ and $T$. In the 
first case, we have two possiblities; $b_6$ $=$ $b_\a$,
$b_\b$ and $b_7$ $=$ $b_\e$. In the former (resp. latter) 
case, we say that the block has type $\B$ (resp. $\C$) 
and say that a block of type $Y_{4,6}^{\cdot,a}$ has 
type $\B'$ (resp. $\C'$) if it represents a mirror image 
of the diagram for $\B$ (resp. $\C$) with $T$. In the 
second case, we also have two possiblities; $b_6$ $=$ 
$b_\c$ or $b_7$ $=$ $b_\c$, since the boundary curve of 
its $f_4$ is surrounded by the boundary curve of its 
$f_6$ on the diagram. In the fomer case, we say that the 
block has type $\D$ and define type $\D'$ as above. The 
latter case contradicts the primeness of $\t{L}$. In the
third case, we have two possibilities; $b_7$ $=$ $b_\d$,
$b_\e$ and $b_6$ $=$ $b_\a$. The former case contradicts 
the reducedness and the latter case contradicts the 
minimality of $\t{L}$ (diagram V).
%%%%%%%%%%%%%%%%%%%%%%%%%%%%%%%%%%%
%
%%%%%%%%%%%%%%%%
\par\vspace{3mm}
%%%%%%%%%%%%%%%%
%
%%%%%%%%%%%%%%%%%%%%%%%%%%%%%%%%%%%
Take a look at $Y_{6,6}^{a,a}$ and put names $v_1$, $\ldots$, 
$v_9$ as shown in Figure \ref{fig:y46y66z86}. Call the face 
with $v_1$ (resp. $v_9$) a face $f$ (resp. $f'$). Let $T_1$
and $T_2$ be two blocks of type $T$. Let $f$ (resp. $f'$) be 
$A$-adjacent to $T_1$ (resp. $T_2$). Here we assume that $T_1$ 
(resp. $T_2$) has vertices $v_\a$, $v_\b$, $v_\c$, $v_\d$, and 
$v_\e$ (resp. $v_{\a'}$, $v_{\b'}$, $v_{\c'}$, $v_{\d'}$, and 
$v_{\e'}$). From the symmetricity, we may assume that the boundary 
curve of $f$ passes the rightside of the dealternator on the 
diagram and then, $f$ is $A$-adjacent to $T_1$ at $v_3$. Replace 
the boundary cycles of $f$ and $T_1$ on the diagram. Now we have 
two possibilities to replace the boundary cycle of $f'$; $b_6$ $=$ 
$b_\c$ or $b_8$ $=$ $b_\c$ from the almost alternating property. 
In the first case, we have that $b_7$ $=$ $b_{\b'}$. If $b_{\a'}$ 
$=$ $b_\d$ or $b_{\b'}$ $=$ $b_1$, then it contradicts the 
reducedness. Therefore, we have that $b_{\a'}$ $=$ $b_1$ or 
$b_{\b'}$ $=$ $b_\d$. In the former (resp. latter) case, we say 
that the block has type $\E$ (resp. $\F$) and define type $\E'$ 
and type $\F'$ as before. In the second case, we have that $b_7$ 
$=$ $b_{\d'}$ and $b_9$ $=$ $b_{\b'}$. If $b_{\a'}$ $=$ $b_\d$, 
then it contradicts the minimality of $\t{L}$ (diagram IV).
If $b_{\b'}$ $=$ $b_\d$, then it contradicts the primeness. Thus, 
we may assume that $b_{\a'}$ $=$ $b_1$ and that $b_{\d'}$ $=$ $b_\b$ 
or $b_{\e'}$ $=$ $b_\b$. If $b_{\e'}$ $=$ $b_\b$, then it 
contradicts the minimality, again. Therefore, we have that 
$b_{\e'}$ $=$ $b_\a$ and $b_{\d'}$ $=$ $b_\b$, and then we say 
that the block has type $\G$ and define type $\G'$ as above.   
%%%%%%%%%%%%%%%%%%%%%%%%%%%%%%%%%%%
%
%%%%%%%%%%%%%%%%
\par\vspace{3mm}
%%%%%%%%%%%%%%%%
%
%%%%%%%%%%%%%%%%%%%%%%%%%%%%%%%%%%%       
At last, take a look at $Z_{8,6}^{a,a}$ and put names $v_1$, 
$\ldots$, $v_{10}$ as shown in Figure \ref{fig:y46y66z86}. 
Let its $f_6$ and its $f_8$ be $A$-adjacent to, a block of
type $T$, $T_1$ and $T_2$, respectively. Replace the boundary 
cycle of its $f_8$ on the diagram assuming that we have already 
replaced the boundary cycles of the $f_6$ and $T_1$. First, assume 
that the boundary curve of the $f_6$ passes the rightside of 
the dealternator on the diagram, and then it is $A$-adjacent 
to $T_1$ at $v_9$. Then, its $f_8$ must be $A$-adjacent to $T_2$ 
at $v_4$ from the almost alternating property and the minimality
of $\t{L}$ (diagram IV). Then, we have that $b_3$ $=$ $b_{\d'}$ 
and $b_5$ $=$ $b_{\b'}$. If $b_{\a'}$ $=$ $b_\d$ or $b_{\e'}$ 
$=$ $b_\b$, then it contradicts the minimality of $\t{L}$
(diagram IV). Thus, we have that $b_{\b'}$ $=$ $b_\d$ and 
$b_{\d'}$ $=$ $b_\b$. Then, we say that the block has type $\H$. 
Second, assume that the boundary curve of the $f_6$ passes the 
leftside of the dealternator on the diagram, and then we may 
assume that it is $A$-adjacent to $T_1$ at $v_{10}$ from the 
symmetricity of the block. Then, its $f_8$ is $A$-adjacent to 
$T_2$ at $v_3$ or $v_5$ from the almost alternating property. 
In the first case, it contradicts the reducedness if $b_4$ $=$ 
$b_\e$ and the primeness if $b_5$ $=$ $b_\e$. It also contradicts 
the the minimality of $\t{L}$ (diagram IV) if $b_6$ $=$ $b_\d$ or 
$b_6$ $=$ $b_\e$. Therefore, we have that $b_4$ $=$ $b_\d$ or 
$b_5$ $=$ $b_\d$. In the former (resp. latter) case, we call the 
block has type $\I$ (resp. $\J$). In the second case, we have that 
$b_6$ $=$ $b_\d$ or $b_6$ $=$ $b_\e$. Both cases contradict the 
minimality of $\t{L}$ (diagram IV).
%%%%%%%%%%%%%%%%%%%%%%%%%%%%%%%%%%%
%
%%%%%%%%%%%%%%%%%%%%%%%%%%%%%%%%%%%
\par\vspace{3mm}\begin{figure}[htbp]
\caption{$Y_{4,6}$, $Y_{6,6}$, and $Z_{8,6}$}
\label{fig:y46y66z86}\end{figure}\par\vspace{3mm}
%%%%%%%%%%%%%%%%%%%%%%%%%%%%%%%%%%%
%
%%%%%%%%%%%%%%%%%%%%%%%%%%%%%%%%%%%
\begin{figure}[htbp]\caption{}
\label{fig:AtoJ}\end{figure}\par\vspace{3mm}
%%%%%%%%%%%%%%%%%%%%%%%%%%%%%%%%%%%
%
%%%%%%%%%%%%%%%%%%%%%%%%%%%%%%%%%%%%%%%%%%% CLAIM ClmCoexst %%%%%
\begin{Claim}\label{ClmCoexst}                                   
Let $L$ $=$ $\{\A$, $\B$, $\C\}$, $L'$ $=$ $\{\B'$, $\C'\}$, 
    $M$ $=$ $\{\G$, $\H\}$, $M'$ $=$ $\{\G'\}$, 
    $N$ $=$ $\{\D$, $\E$, $\F$, $\I$, $\J\}$, and 
    $N'$ $=$ $\{\D'$, $\E'$, $\F'$, $\I'$, $\J'\}$.            
Then, we have the following. 
%%%%%%%%%%%%%%%%%%%%%%%%%%%%%%%%%%%      
\begin{itemize}   
\item[(i)]  A block of any type of $L$ $\cup$ $L'$ $\cup$ $M$ $\cup$
            $M'$ does not coexist in graph $G$ with any other 
            block of a different type. Blocks of the same type 
            of $L$ $\cup$ $L'$ $\cup$ $M$ $\cup$ $M'$ can coexist 
            in graph $G$. Then, their boundary curves of their top 
            $($resp. bottom$)$ faces pass the same four bubbles.
\item[(ii)] Any block of $N$ does not coexist in the graph $G$ 
            with any blocks of $N'$. Blocks of $N$ $($resp. 
            $N'$$)$ may coexist in graph $G$. Then, the boundary 
            curves below $($resp. above$)$ the dealternator
            pass the same five bubbles.   
\end{itemize}\end{Claim}       
%%%%%%%%%%%%%%%%%%%%%%%%%%%%%%%%%%%%%%%%%%%%%%%%% END CLAIM %%%%%
%
%%%%%%%%%%%%%%%%%%%%%%%%%%%%%%%%%%%%%%%%%%%%%%%%%%%%% PROOF %%%%%
\begin{proof}
%%%%%%%%%%%%%%%%%%%%%%%%%%%%%%%%%%%%%%%%%%%%%%% L & M and N %%%%%
(i) We say a subgraph of $G$ a subblock of type $P$, $Q$, $R$, 
and $S$ (resp. $P'$, $Q'$, $R'$, and $S'$) if its boundary curve
constracts a diagram $P$, $Q$, $R$, and $S$ (resp. the mirror 
image $P'$, $Q'$, $R'$, and $S'$) of Figure \ref{fig:PQRS}, 
respectively. It is easy to see that it contradicts the mimimality 
of $\t{L}$ (diagram IV) if graph $G$ has $S$ (resp. $S'$) and 
$P'$, $Q'$ and $R'$ (resp. $P$, $Q$ and $R$). Here note that any 
block of a type of $L$ $\cup$ $L'$ $\cup$ $M$ $\cup$ $M'$ consists 
of one of $P$, $Q$ and $R$ and one of $P'$, $Q'$ and $R'$ (for 
instance, a block of type $\G$ consists of a subblock of type $P$ 
and a subblock of type $R'$). In addition, any block of $N$ (resp. 
$N'$) contains a subblock of type $S'$ (resp. $S$). Therefore, any 
block of $L$ $\cup$ $L'$ $\cup$ $M$ $\cup$ $M'$ and any block of 
$N$ $\cup$ $N'$ cannot coexist in graph $G$.
%%%%%%%%%%%%%%%%%%%%%%%%%%%%%%%%%%%   
%
%%%%%%%%%%%%%%%%%%%%%%%%%%%%%%%%%%%%%%% L & M ( P outside ) %%%%%
From the primeness, the type of a block whose boundary curve can 
exist inside of the boundary curve of $P$ is only $P$ among $P$, 
$Q$, and $R$, and then their boundary curves pass the same 
bubbles. 
%%%%%%%%%%%%%%%%%%%%%%%%%%%%%%%%%%%   
%
%%%%%%%%%%%%%%%%%%%%%%%%%%%%%%%%%%%%%%% L & M ( Q outside ) %%%%%
Next, assume that there is the boundary curve 
$b_d$$b_\c$$b_\d$$b_\e$ of a face of degree $4$ inside of a
boundary curve $b_d$$b_{\c'}$$b_{\d'}$$b_{\e'}$ of a subblock of
type $Q$. Then we have that 
$b_\d$ $=$ $b_{\d'}$ and $b_\e$ $=$ $b_{\e'}$,
$b_\d$ $=$ $b_1$ and $b_\e$ $=$ $b_{\e'}$,
$b_\d$ $=$ $b_2$, or $b_\e$ $=$ $b_2$.
The last three cases contradicts the primeness. Therefore
consider the first case. If the face is of a subblock of type $P$,
then it contradicts primeness. If the face is of a subblock of type
$Q$, then their boundary curves pass the same $6$ bubbles from
the primeness. If the face is of a subblock of type $R$, then we
cannot connect bubbles $b_4$ and $b_\d$ with an arc for $R$ (see 
Figure \ref{fig:PQRS}).
%%%%%%%%%%%%%%%%%%%%%%%%%%%%%%%%%%%   
%
%%%%%%%%%%%%%%%%%%%%%%%%%%%%%%%%%%%%%%% L & M ( R outside ) %%%%%
Now assume that there is the boundary curve 
$b_d$$b_\c$$b_\d$$b_\e$ of a face of degree $4$ inside of the
boundary curve $b_d$$b_{\c''}$$b_{\d''}$$b_{\e''}$ of a subblock 
of type $R$. Then we have that $b_\d$ $=$ $b_{\d''}$, 
$b_\e$ $=$ $b_{\d''}$, $b_\d$ $=$ $b_4$, or $b_\e$ $=$ $b_4$.
The last two cases contradicts the primeness. Consider the first
case. Then, we also obtain that $b_\e$ $=$ $b_{\e''}$ from the
primeness. If the face is of a subblock of type $P$ or $Q$, then
it contradicts the primeness. If the face is of a subblock of 
type $R$, then their boundary curves pass the same $6$ bubbles 
also from the primeness. Next, consider the second case. Note 
that we are now considering the coexistence of blocks of $L$ 
$\cup$ $L'$ $\cup$ $M$ $\cup$ $M'$. Therefore, we have a 
subblock of type $P'$, $Q'$, or $R'$. Then it contradicts the 
minimality of $\t{L}$ (diagram IV).
%%%%%%%%%%%%%%%%%%%%%%%%%%%%%%%%%%% 
%
%%%%%%%%%%%%%%%%%%%%%%%%%%%%%%%%%%%  
Now we need to show that a block of type $\C$ and a block of 
type $\G$ (or $\C'$ and $\G'$) do not coexist in graph $G$. 
If graph $G$ has $\C$, then the boundary curve of the top and 
bottom face of any block of type $T$ must pass the same $4$ 
bubbles as the boundary curve of the top and bottom face of 
the block of type $\C$, respectively. However, $\G$ has two 
top faces whose boundary curves do not pass the same bubbles. 
It is a contradiction. (ii) If any of $N$ and any of $N'$ 
coexist in graph $G$, then we have $S$ and $S'$ on the diagram, 
which contradicts the reducedness. It is easy to see the last 
part following the previous case.
\end{proof}           
%%%%%%%%%%%%%%%%%%%%%%%%%%%%%%%%%%%%%%%%%%%%%%%%% END PROOF %%%%%
%
%%%%%%%%%%%%%%%%%%%%%%%%%%%%%%%%%%% 
\par\vspace{3mm}\begin{figure}[htbp]\caption{}
\label{fig:PQRS}\end{figure}\par\vspace{3mm}
%%%%%%%%%%%%%%%%%%%%%%%%%%%%%%%%%%% 
%
%%%%%%%%%%%%%%%%%%%%%%%%%%%%%%%%%%%
We devide Case 1 into the following $6$ subcases; case 1-$\A$, 
$\B$, $\C$, $\G$, $\H$, and $N$ according that there is a block 
of $\A$, $\B^*$, $\C^*$, $\G^*$, $\H^*$, and $N$, respectively. 
In each case, we look at the block which is $B_*$-adjacent to 
a block of type $T$ with a negative weight. Then, We have the 
following.
%%%%%%%%%%%%%%%%%%%%%%%%%%%%%%%%%%%
%
%%%%%%%%%%%%%%%%%%%%%%%%%%%%%%%%%%%%%%%%%%% CLAIM Clm1Small %%%%%
\begin{Claim}\label{Clm1Small}
%%%%%%%%%%%%%%%%%%%%%%%%%%%%%%%%%%%
\begin{itemize}
\item[(i)]    No face of degree $4$ can be $B_b$- $($or $B_t$-$)$ 
              adjacent to a block of type $T$.
\item[(ii)]   No face of degree $6$ can be $B_b$- and 
              $B_t$-adjacent to a block or blocks of type $T$.
\item[(iii)]  No face of degree $8$ can be $A$-, $B_b$-, and 
              $B_t$-adjacent to blocks of type $T$.
\end{itemize}\end{Claim}
%%%%%%%%%%%%%%%%%%%%%%%%%%%%%%%%%%%%%%%%%%%%%%%%% END CLAIM %%%%%
%
%%%%%%%%%%%%%%%%%%%%%%%%%%%%%%%%%%%%%%%%%%%%%%%%%%%%% PROOF %%%%%
\begin{proof}
(i) Assume that we have a face of degree $4$ with boundary 
cycle $v_\a$$v_\b$$v_1$$v_2$. Considering the length of the
cycle, we have that $v_1$ $=$ $v_d$ and $b_2$ $=$ $b_\d$ or
$b_\e$. The former case contradicts the reducedness and the 
latter case contradicts the mimimality of $\t{L}$ (diagram V).
(ii) Assume that we have a face of degree $6$ which is $B_b$- 
and $B_t$-adjacent to a block of type $T$ and its boundary 
cycle is $v_\a$$v_\b$$v_1$$v_2$$v_3$$v_4$ (it can be similarly 
shown the case that the face is $B_b$- and $B_t$-adjacent to 
blocks of type $T$). Considering the length of the cycle and 
the almost alternating property, we have that $v_1$ $=$ $v_d$, 
$b_2$ $=$ $b_\d$, and $b_3$ $=$ $b_\e$. Then, it contradicts 
the reducedness.
(iii) We show only the case that we have a face of degree $8$ 
which is $A$-adjacent to a block of type $T$ and $B_b$-, and 
$B_t$-adjacent to another block of type $T$. And let 
$v_\a$$v_\b$$v_1$$v_2$$v_3$$v_4$$v_5$$v_6$ be the boundary
cycle of the face. Considering the length of the cycle and
the almost alternating property, we have that the face is 
$B_t$-adjacent to the block at $v_2$$v_3$ or $v_4$$v_5$.
The former case contradicts the reducedness. In the latter 
case, the face must be $A$-adjacent to the block at $v_6$,
and then it contradicts the reducedness, again.
\end{proof}
%%%%%%%%%%%%%%%%%%%%%%%%%%%%%%%%%%%%%%%%%%%%%%%%% END PROOF %%%%%
%
%%%%%%%%%%%%%%%%
\par\vspace{3mm}
%%%%%%%%%%%%%%%%
%
%
%%%%%%%%%%%%%%%%%%%%%%%%
\noindent{\bf Case 1-$\A$.}
%%%%%%%%%%%%%%%%%%%%%%%%
%
%%%%%%%%%%%%%%%%%%%%%%%%%%%%%%%%%%%
Take a block of type $\A$ and a block of type $T$ which are 
$A$-adjacent to each other and put them together. We call it 
a block of type $T_{\A}$ or simply $T_{\A}$, and so $w(T_{\A})$ 
$=$ $-2$. Note that we do not have any of $\{L$ $\cup$ $L'$ 
$\cup$ $M$ $\cup$ $M'$ $\cup$ $N$ $\cup$ $N'\}$ $-$ $\{\A\}$
from Claim \ref{ClmCoexst}. Take a look at blocks which are 
$B_*$-adjacent to blocks of type $T_{\A}$. We define, for 
instance, $Z_{i, j}^{a*, bt}$ as a block whose $f_i$ is $A$- 
and $B_*$-adjacent to blocks of type $T_{\A}$ and whose $f_j$ 
is $B_b$- and $B_t$-adjacent to blocks of type $T_{\A}$. Also, 
we use $Z_i^a$ for a face $f_i$ which is adjacent to two faces 
of degree $2$ and is $A$-adjacent to $T_{\A}$. 
%%%%%%%%%%%%%%%%%%%%%%%%%%%%%%%%%%%
%
%%%%%%%%%%%%%%%%
\par\vspace{3mm}
%%%%%%%%%%%%%%%%
%
%%%%%%%%%%%%%%%%%%%%%%%%%%%%%%%%%%%%%%%%%%% CLAIM Clm1aSmall %%%%
\begin{Claim}\label{Clm1aSmall} Graph $G$ does not have any
block of type $X_6^{a*}$, $X_{10}^{abt}$, $Z_4^{\cdot}$, $Z_6^a$, 
$Z_8^{a*}$, $Y_{8, j}^{bt, a}$, or $Z_{6, j}^{*, a}$. 
\end{Claim}
%%%%%%%%%%%%%%%%%%%%%%%%%%%%%%%%%%%%%%%%%%%%%%%%% END CLAIM %%%%%
%
%%%%%%%%%%%%%%%%%%%%%%%%%%%%%%%%%%%%%%%%%%%%%%%%%%%%% PROOF %%%%%
\begin{proof}
Let $v_\a$$v_\b$$v_1$$v_2$$v_3$$v_4$ be the boundary cycle of 
$X_6^{ab}$. Note that the boundary curve must pass the 
dealternator and curve $b_d$$b_\c$$b_\d$$b_\e$. Thus, $X_6^{ab}$
must be $A$-adjacent to a block of type $T$ at $v_4$ from the
almost alternating property. However then, it contradicts the
minimality of $\t{L}$ (diagram IV).
We can similarly show that there does not exist $X_{10}^{abt}$.
%%%%%%%%%%%%%%%%%%%%%%%%%%%%%%%%%%%
%
%%%%%%%%%%%%%%%%
\par\vspace{3mm}
%%%%%%%%%%%%%%%%
%
%%%%%%%%%%%%%%%%%%%%%%%%%%%%%%%%%%%
Let $v_d$$v_1$$v_2$$v_3$ be the boundary cycle of $Z_4^{\cdot}$
and let $b_1$ (resp. $b_3$) be inside of $b_d$$b_\a$$b_\b$$b_\c$
(resp. $b_d$$b_\c$$b_\d$$b_\e$). Since the length of curve
$b_d$$b_1$$b_2$$b_3$ is $4$ and it must pass 
$b_d$$b_\a$$b_\b$$b_\c$ and $b_d$$b_\c$$b_\d$$b_\e$, we have 
that $b_1$ $=$ $b_\a$, $b_\b$ or that $b_3$ $=$ $b_\d$, $b_\e$.
However if $b_1$ $=$ $b_\a$ or $b_3$ $=$ 
$b_\e$, then it contradicts the primeness. And if $b_1$ $=$ 
$b_\b$ or $b_3$ $=$ $b_\d$, then it contradicts the minimality 
of $\t{L}$ (diagram IV). In the case of 
$Z_6^a$, let $v_d$$v_1$$v_2$$v_3$$v_4$$v_5$ be its boundary 
cycle. Following the previous argument, we have that $b_1$ $\ne$
$b_\a$, $b_\b$ and $b_5$ $\ne$ $b_\d$, $b_\e$. From the length 
of the cycle, $Z_6^a$ must be $A$-adjacent to a block of type
$T$ at $v_3$. Then, we have that $b_2$ $=$ $b_\a$ or $b_\b$.
The first case contradicts the minimality of $\t{L}$
(diagram IV) and the second case contradicts the primeness.
We can similarly show that there does not exist $Z_8^{a*}$.
%%%%%%%%%%%%%%%%%%%%%%%%%%%%%%%%%%%
%
%%%%%%%%%%%%%%%%
\par\vspace{3mm}
%%%%%%%%%%%%%%%%
%
%%%%%%%%%%%%%%%%%%%%%%%%%%%%%%%%%%%
Let $v_\a$$v_\b$$v_1$$v_2$$v_3$$v_4$$v_5$$v_6$ be the boundary 
cycle of $Y_8^{bt}$. It is $B_t$-adjacent to $T$ at $v_2$$v_3$
or $v_4$$v_5$ from the almost alternating property. The former
case contradicts the reducedness. In the latter case, we have
that $b_3$ $=$ $b_d$ and $b_6$ $\ne$ $b_\c$ from the reducedness
and the minimality of $\t{L}$ (diagram IV). Therefore, $b_\c$ 
and the boundary curve of $f_j$ of $Y_{8,j}^{bt,a}$ is in one 
and in the other of the two regions of 
$S^2$ $-$ $b_\a$$b_\b$$b_1$$b_2$$b_3$$b_4$$b_5$$b_6$, 
respectively. Therefore, $f_j$ cannot be $A$-adjacent to $T$. 
We can similarly show that there does not exist $Z_{6,j}^{*,a}$.
\end{proof}
%%%%%%%%%%%%%%%%%%%%%%%%%%%%%%%%%%%%%%%%%%%%%%%%% END PROOF %%%%%
%
%%%%%%%%%%%%%%%%
\par\vspace{3mm}
%%%%%%%%%%%%%%%%
%
%%%%%%%%%%%%%%%%%%%%%%%%%%%%%%%%%%%%%%%%%%%% CLAIM Clm1aLeft %%%%
\begin{Claim}\label{Clm1aLeft} The boundary curve of the face 
$($$\ne$ $f_2$$)$ of any of $X_8^{a*}$, $X_8^{bt}$, $Y_6^a$, and 
$Z_6^*$ passes the leftside of the dealternator on the diagram. 
\end{Claim}
%%%%%%%%%%%%%%%%%%%%%%%%%%%%%%%%%%%%%%%%%%%%%%%%% END CLAIM %%%%%
%
%%%%%%%%%%%%%%%%%%%%%%%%%%%%%%%%%%%%%%%%%%%%%%%%%%%%% PROOF %%%%%
\begin{proof}
We show the proof only for $X_8^{ab}$ and $Y_6^a$. Let 
$v_\a$$v_\b$$v_1$$v_2$$v_3$$v_4$$v_5$$v_6$ be the boundary 
cycle of $X_8^{ab}$. Note that the boundary curve must pass 
the dealternator and curve $b_d$$b_\c$$b_\d$$b_\e$. Thus, 
$X_8^{ab}$ must be $A$-adjacent to a block of type $T$ at 
$v_4$ or $v_6$ from the almost alternating property. Then
we have the former case and that $v_2$ $=$ $v_d$
from the minimality of $\t{L}$ (diagram IV), which completes 
the proof. Moreover, we have that $b_3$ $=$ $b_\d$ from the 
primeness and then we obtain the diagram shown in Figure 
\ref{fig:x8ay6a}.
%%%%%%%%%%%%%%%%%%%%%%%%%%%%%%%%%%%
%
%%%%%%%%%%%%%%%%
\par\vspace{3mm}
%%%%%%%%%%%%%%%%
%
%%%%%%%%%%%%%%%%%%%%%%%%%%%%%%%%%%%
Let $v_{\c'}$$v_1$$v_2$$v_3$$v_4$$v_5$ be the boundary cycle 
of $Y_6^a$. Since the length of the boundary curve is $6$ and 
it must pass the boundary curves of the top and bottom faces 
of the block of type $T$ of $T_\A$, we have that $b_d$ $=$ 
$b_2$, $b_3$, or $b_4$. If $b_d$ $=$ $b_3$, then we have that 
$b_1$ $=$ $b_\a$, $b_\b$ or that $b_2$ $=$ $b_\a$, $b_\b$. 
The first and fourth cases contradict the mimimality of $\t{L}$ 
(diagram IV) and the second case contradicts the primeness. 
Therefore we have that $b_2$ $=$ $b_\a$ and then, we similarly 
obtain that $b_4$ $=$ $b_\e$. However, then it contradicts the 
primeness, since $Y_6^a$ is adjacent to a face of degree $2$ at 
$v_2$$v_3$ or at $v_3$$v_4$. This completes the proof. Moreover, 
if $b_d$ $=$ $b_2$, then $b_1$ $=$ $b_\a$ or $b_\b$. The former 
case contradicts the mimimality of $\t{L}$ (diagram IV). In the 
latter case, $Y_6^a$ must be adjacent to a face of degree $2$ 
at $v_2$$v_3$ from the minimality of $\t{L}$ (diagram IV). Then, 
we have that $b_3$ $=$ $b_\d$, $b_\e$, that $b_4$ $=$ $b_\d$, 
$b_\e$, or that $b_5$ $=$ $b_\d$, $b_\e$. All the five cases 
but the third contradict the primeness or the minimality of 
$\t{L}$ (diagram IV). Therefore, we obtain the diagram shown 
in Figure \ref{fig:x8ay6a}. In the case that $b_4$ $=$ $b_d$, 
we obtain a mirror image of the diagram.
\end{proof}
%%%%%%%%%%%%%%%%%%%%%%%%%%%%%%%%%%%%%%%%%%%%%%%%% END PROOF %%%%%
%
%%%%%%%%%%%%%%%%%%%%%%%%%%%%%%%%%%% 
\par\vspace{3mm}\begin{figure}[htbp]
\caption{$X_8^{a*}$ and $Y_6^a$ with $T_\A$ and $T'$}
\label{fig:x8ay6a}\end{figure}\par\vspace{3mm}
%%%%%%%%%%%%%%%%%%%%%%%%%%%%%%%%%%% 
%
%%%%%%%%%%%%%%%%%%%%%%%%%%%%%%%%%%% 
Claim \ref{Clm1aLeft} says that there does not exist any 
$X_{i,j}$ or $Y_{i,j}$ such that each of its faces $f_i$ and 
$f_j$ is $X_8^{a*}$, $X_8^{bt}$, $Y_6^a$, or $Z_6^*$. Moreover, 
we have the following.
%%%%%%%%%%%%%%%%%%%%%%%%%%%%%%%%%%% 
%
%%%%%%%%%%%%%%%%
\par\vspace{3mm}
%%%%%%%%%%%%%%%%
%
%%%%%%%%%%%%%%%%%%%%%%%%%%%%%%%%%%%%%%%%% CLAIM Clm1aSmlPair %%%%
\begin{Claim}\label{Clm1aSmlPair} Graph $G$ does not have 
any block of type $Y_{4, 8}^{\cdot, a*}$, $Y_{6, 6}^{a, *}$, 
or $Y_{6, 8}^{*, a*}$.
\end{Claim}
%%%%%%%%%%%%%%%%%%%%%%%%%%%%%%%%%%%%%%%%%%%%%%%%% END CLAIM %%%%%
%
%%%%%%%%%%%%%%%%%%%%%%%%%%%%%%%%%%%%%%%%%%%%%%%%%%%%% PROOF %%%%%
\begin{proof}
Take a look at the diagram of $X_8^{a*}$ with $T_\A$ and $T'$
(Figure \ref{fig:x8ay6a}). From the mimimality of $\t{L}$
(diagram IV), $Y_8^{a*}$ must be adjacent to a face of degree 
$2$ at $b_1$$b_2$. Then, it is easy to see that there does not 
exist $Y_{4,8}^{\cdot, a*}$, since the length of the boundary 
curve of $Y_4^{\cdot}$ is $4$ and it must pass 
$b_d$$b_\a$$b_\b$$b_\c$, $b_d$$b_{\a'}$$b_{\b'}$$b_{\c'}$, and 
$b_d$$b_\c$$b_\d$$b_\e$. We can similarly show that there does 
not exist $Y_{6,6}^{a,*}$ or $Y_{6,8}^{*,a*}$.
\end{proof}
%%%%%%%%%%%%%%%%%%%%%%%%%%%%%%%%%%%%%%%%%%%%%%%%% END PROOF %%%%%
%
%%%%%%%%%%%%%%%%
\par\vspace{3mm}
%%%%%%%%%%%%%%%%
%
%%%%%%%%%%%%%%%%%%%%%%%%%%%%%%%%%%%
We are now looking at blocks which are 
$B_*$-adjacent to blocks of type $T_{\A}$. Each face ($\ne$ 
$f_2$) of every such a block can be $B_b$- (resp. $B_t$-) 
adjacent to at most one $T_{\A}$ at most once from Claim 
\ref{ClmCoexst} (i). In addition, the face might be $A$-adjacent 
to a block of type $T$ as well, and then note that we have 
discharged weight $2$ of the face to the block of type $T$.
Then we have $68$ types of blocks which are $B_*$-adjacent to 
blocks of type $T_{\A}$; $X_i^p$, $Y_{i,j}^{q,r}$, and 
$Z_{i,j}^{q,r}$ with $\{q$, $r\}$ $=$ $\{\cdot$, $*$, $bt$, $a$, 
$a*$, $abt\}$ and $p$ and one of $q$ and $r$ are of $\{*$, $bt$,
$a*$, $abt\}$. Consider the sum of weights of such a block and 
blocks of type $T$ or type $T_{\A}$ which it is $A$- or 
$B_*$-adjacent to. Then the type of a block such that the sum
is negative is 
$Y_{4,6}^{\cdot,*}$, $Y_{6,4}^{*,\cdot}$,
$Z_{6,6}^{\cdot,*}$, $Y_{6,4}^{*,\cdot}$,
$Y_{6,6}^{*,*}$, $Z_{8,6}^{*,*}$, $Z_{10,8}^{bt,bt}$,
$Y_{4,8}^{\cdot,bt}$, $Y_{8,4}^{bt,\cdot}$, $Z_{8,6}^{bt,\cdot}$,
$Y_{6,8}^{*,bt}$, $Y_{8,6}^{bt,*}$, 
$Z_{8,8}^{*,bt}$, $Z_{8,8}^{bt,*}$, or 
$Z_{10,6}^{bt,*}$ from Claim \ref{ClmAadj}, \ref{Clm1Small}, 
\ref{Clm1aSmall}, \ref{Clm1aLeft}, and \ref{Clm1aSmlPair}   
and then the sum is $-2$. Next we consider such blocks.
%%%%%%%%%%%%%%%%%%%%%%%%%%%%%%%%%%%
%
%%%%%%%%%%%%%%%%
\par\vspace{3mm}
%%%%%%%%%%%%%%%%
%
%%%%%%%%%%%%%%%%%%%%%%%%%%%%%%%%%%% 
In each case of $Y_{6,6}^{b,b}$ and $Y_{6,6}^{t,t}$, choose one 
of two blocks of type $T_{\A}$ and discharge its weight $2$ to 
the block. In the case of $Z_{8,6}^{*,*}$, discharge the weight 
$2$ to the block of type $T_{\A}$ which its $f_8$ is 
$B_*$-adjacent to. Now, in each of the above $3$ cases and the 
first $4$ of $16$ cases, we have the situation that a block with 
its weight $0$ is $B_*$-adjacent to $T_{\A}$ with its weight $-2$. 
We say that such blocks are type I. In each of the last $5$ of 
$16$ cases, discharge $2$ out of its weight $4$ to the block of 
type $T_{\A}$ which it is $B_*$-adjacent to. In the case of 
$Z_{10,8}^{bt,bt}$, discharge $2$ out of its weight $6$ to each 
of the two blocks of type $T_{\A}$ which its $f_{10}$ is $B_b$- 
and $B_t$-adjacent to. Now, in each of the above $6$ cases and 
the cases of $Y_{4,8}^{\cdot,bt}$, $Y_{8,4}^{bt,\cdot}$, 
$Y_{6,6}^{b,t}$, $Y_{6,6}^{t,b}$, and $Z_{8,6}^{bt,\cdot}$, we 
have that a block with its weight $2$ is $B_b$- and 
$B_t$-adjacent to two blocks of type $T_{\A}$ with each weight 
$-2$ (if the two blocks of type $T_{\A}$ are the same, we can 
discharge the weight $2$ to the block of type $T_{\A}$ and 
make its weight non-negative. Therefore, we do not consider  
such a case). We say that such blocks are type II.  
%%%%%%%%%%%%%%%%%%%%%%%%%%%%%%%%%%%
%
%%%%%%%%%%%%%%%%
\par\vspace{3mm}
%%%%%%%%%%%%%%%%
%
%%%%%%%%%%%%%%%%%%%%%%%%%%%%%%%%%%%
Then, we can constract paths by regarding blocks of type $T$ with 
negative weights and blocks of type I and type II as edges and 
vertices, respectively. Here, note that each block of type I is 
$B_*$-adjacent to exactly one $T$ with a negative weight and 
each block of type II is $B_b$- and $B_t$-adjacent to exactly 
two blocks of type $T$ with negative weights. Therefore, for 
each path, if the block corresponding to one of its ends is 
$B_b$-adjacent to $T$ with a negative weight, then the block 
corresponding to the other of its ends is $B_t$-adjacent to $T$ 
with a negative weight. Now we have the following.
%%%%%%%%%%%%%%%%%%%%%%%%%%%%%%%%%%%
%
%%%%%%%%%%%%%%%%
\par\vspace{3mm}
%%%%%%%%%%%%%%%%
%
%%%%%%%%%%%%%%%%%%%%%%%%%%%%%%%%%%%%%%%%%%%%% CLAIM ClmArcs %%%%%
\begin{Claim}\label{ClmArcs}                                  
Assume that diagram $D$ contains at least one from each of 
$\{P$, $Q$, $R\}$. Let $\x$ $($resp. 
$\x'$$)$ be an arc $b_1$$b_2$$b_3$$b_4$ such that $b_1$ $=$ 
$b_\b$, $b_2$ $=$ $b_\a$, and $b_4$ $=$ $b_\d$ $($resp. $b_1$ 
$=$ $b_\d$, $b_2$ $=$ $b_\e$, and $b_4$ $=$ $b_\b$$)$. Let $\y$ 
$($resp. $\y'$$)$ be an arc $b_1$$b_2$$b_3$ such that $b_1$ $=$ 
$b_\b$, $b_2$ $=$ $b_\a$, and $b_3$ $=$ $b_\d$ $($resp. $b_1$ 
$=$ $b_\d$, $b_2$ $=$ $b_\e$, and $b_3$ $=$ $b_\b$$)$. Let $\z$ 
$($resp. $\z'$$)$ be an arc $b_1$$b_2$ such that $b_1$ $=$ 
$b_\b$ and $b_2$ $=$ $b_\d$ $($resp. $b_1$ $=$ $b_\d$ and $b_2$ 
$=$ $b_\b$$)$ and it sees the center of $b_\b$ on the same 
$($resp. opposite$)$ side as it does $b_\c$ $($Figure 
\ref{fig:arcs}$)$. Suppose that $D$ contains one of $\{\x$, $\y$, 
$\z\}$. Then, $D$ contains none of $\{\x'$, $\y'$, $\z'\}$.
\end{Claim}                                                      
%%%%%%%%%%%%%%%%%%%%%%%%%%%%%%%%%%%%%%%%%%%%%%%%% END CLAIM %%%%%
%
%%%%%%%%%%%%%%%%%%%%%%%%%%%%%%%%%%%%%%%%%%%%%%%%%%%%% PROOF %%%%%
\begin{proof}
Assume that $D$ contains $\x$. Then, any arc passing $b_\b$ and
$b_\e$ must passes $b_\a$. Thus, $D$ does not contain $\x'$ or 
$\y'$. Since $\z'$ passes $b_\b$ (resp. $b_\d$) seeing it at the 
opposite (resp. same) side as $b_\c$, it must be surrounded by 
(resp. it must surround) the cycle containing $\x$. It is a 
contradiction. Thus $D$ does not contain $\z'$. Assume that $D$ 
contains $\y$ and $\z'$. Then the diagram contains $S'$. It 
contradicts the mimimality of $\t{L}$ (diagram IV), since $D$ 
contains $P$, $Q$, or $R$. We can similarly show other cases.
\end{proof}
%%%%%%%%%%%%%%%%%%%%%%%%%%%%%%%%%%%%%%%%%%%%%%%%% END PROOF %%%%%
%
%%%%%%%%%%%%%%%%%%%%%%%%%%%%%%%%%%% 
\par\vspace{3mm}\begin{figure}[htbp]
\caption{$\x$, $\y$, and $\z$}
\label{fig:arcs}\end{figure}\par\vspace{3mm}
%%%%%%%%%%%%%%%%%%%%%%%%%%%%%%%%%%% 
%
%%%%%%%%%%%%%%%%%%%%%%%%%%%%%%%%%%% 
Then, the following claim says that there does not eixst such
a path.
%%%%%%%%%%%%%%%%%%%%%%%%%%%%%%%%%%% 
%
%%%%%%%%%%%%%%%%%%%%%%%%%%%%%%%%%%%%%%%%%% CLAIM Clm1aTpBtm %%%%%
\begin{Claim}\label{Clm1aTpBtm} Assume that there exist a block
of type $T_{\A}$ and a block of type I and they are $B_b$- 
$($resp. $B_t$-$)$ adjacent to each other. Then, there are no 
other blocks of type $T_{\A}$ and no blocks of type I such that 
they are $B_t$- $($resp. $B_b$-$)$ adjacent to each other.    
\end{Claim}      
%%%%%%%%%%%%%%%%%%%%%%%%%%%%%%%%%%%%%%%%%%%%%%%%% END CLAIM %%%%%
%
%%%%%%%%%%%%%%%%%%%%%%%%%%%%%%%%%%%%%%%%%%%%%%%%%%%%% PROOF %%%%%
\begin{proof}
From Claim \ref{ClmArcs}, it is sufficient to show that 
we have an arc $\x$, $\y$, or $\z$ on the diagram under 
the assumptionthat there exists $Y_{4,6}^{\cdot,b}$, 
$Y_{6,4}^{b,\cdot}$, $Y_{6,6}^{b,b}$, or $Z_6^b$. Let 
$f_6$ be $B_b$-adjacent to $T_{\A}$. Then, its boundary 
curve $b_\a$$b_\b$$b_1$$b_2$$b_3$$b_4$ must pass $b_\d$ 
or $b_\e$. In the former case, we have an arc $\x$ or 
$\y$ from the mimimality of $\t{L}$ (diagram IV). In the
latter case, we have that $b_3$ $=$ $b_d$ and $b_4$ $=$ 
$b_\e$ from the mimimality of $\t{L}$ (diagram IV) and 
the primeness. Moreover, it can be adjacent to at most 
one $f_2$ at $v_2$$v_3$ from the primeness. Thus, we do 
not have $Z_6^b$ in this case. If we have $Y_{4,6}^{\cdot,b}$
or $Y_{6,4}^{b,\cdot}$, then the boundary curve of $f_4$ must 
pass $b_\b$ and $b_\d$, and thus we have $\z$. If we have 
$Y_{6,6}^{b,b}$, then the boundary curve of another $f_6$ 
must pass $b_\b$, $b_\a$, and $b_\d$ from the primeness. 
Then, we have $\x$ or $\y$. \end{proof}
%%%%%%%%%%%%%%%%%%%%%%%%%%%%%%%%%%%%%%%%%%%%%%%%% END PROOF %%%%%
%
%%%%%%%%%%%%%%%%
\par\vspace{3mm}                                                 
%%%%%%%%%%%%%%%%                   
%         
%
%%%%%%%%%%%%%%%%%%%%%%%%%%%
\noindent{\bf Case 1-$\B$.} 
%%%%%%%%%%%%%%%%%%%%%%%%%%%
%
%%%%%%%%%%%%%%%%%%%%%%%%%%%%%%%%%%%
We show only the case that we have a block of type $\B$. The case
that we have a block of type $\B'$ can be shown similarly. Define a 
block of type $T_{\B}$ following Case 1-$\A$ and take a look at 
blocks which are $B_*$-adjacent to $T_{\B}$. It is easy to see 
that the boundary curves of the top (resp. bottom) faces of any
blocks of type $T$ pass the same four bubbles as that of the top 
(resp. bottom) face of $T$ of a block of type $T_{\B}$. This 
induces that the boundary curves of the top (resp. bottom) faces
of $T$ of all blocks of type $T_{\B}$ pass the same four babbles.
Therefore, any face which is $A$-adjacent to $T$ cannot be 
$B_*$-adjacent to $T_{\B}$.
Moreover, any face which is $B_*$-adjacent to $T_{\B}$ must be
adjacent to at least one $f_2$, since its boundary curve is outside
of the boundary curve of the face of degree $6$ of $T_{\B}$.
Therefore, we need to take a look at the blocks of type
$Y_{i,j}^{p,q}$ or $Z_{i,j}^{p,q}$, where $p$, $q$ $\in$ 
$\{$ $\cdot$, $a$, $*$, $bt\}$ and one of them is of $\{*$, $bt\}$. 
%%%%%%%%%%%%%%%%%%%%%%%%%%%%%%%%%%%
%
%%%%%%%%%%%%%%%%
\par\vspace{3mm}
%%%%%%%%%%%%%%%%
%
%%%%%%%%%%%%%%%%%%%%%%%%%%%%%%%%%%%
From Claim \ref{ClmAadj} and Claim \ref{Clm1Small}, we know
that $G$ does not have $X_4^*$, $X_6^{bt}$, or $Y_4^a$.
We can similarly show that $G$ does not have $Z_4^{\cdot}$ 
or $Z_6^a$. Moreover, by showing that the boundary curve of 
$Z_6^*$ and $Z_8^{bt}$ passes the leftside of the dealternator 
on the diagram, respectively, we also obtain that $G$ 
does not have any block of type $Z_{6,6}^{*,*}$, 
$Z_{8,6}^{bt,*}$, or $Z_{8,8}^{bt,bt}$. Then the type of a block 
such that the sum of its weight and weights of blocks of type 
$T_{\B}$ which it is $B_*$-adjacent to is negative is 
%%%%%%%%%%
$Y_{6,4}^{*,\cdot}$, $Y_{4,6}^{\cdot,*}$, 
$Z_{6,6}^{*,\cdot}$, $Z_{6,6}^{\cdot,*}$,
$Y_{6,6}^{*,a}$, $Y_{6,6}^{a,*}$, $Z_{8,6}^{a,*}$, 
%%%%%%%%%%
$Y_{8,4}^{bt,\cdot}$, $Y_{4,8}^{\cdot,bt}$,
$Y_{8,6}^{bt,a}$, $Y_{6,8}^{a,bt}$, $Z_{8,6}^{bt,\cdot}$, 
$Z_{8,8}^{bt,a}$, $Z_{8,8}^{a,bt}$,
%%%%%%%%%%
$Y_{6,6}^{*,*}$, $Z_{8,6}^{*,*}$, $Z_{10,6}^{bt,*}$, 
$Y_{8,8}^{bt,bt}$, $Z_{10,8}^{bt,bt}$,
$Y_{8,6}^{bt,*}$, $Y_{6,8}^{*,bt}$,  
$Z_{8,8}^{bt,*}$, or $Z_{8,8}^{*,bt}$,
and then the sum is $-2$.
%%%%%%%%%%%%%%%%%%%%%%%%%%%%%%%%%%%
%
%%%%%%%%%%%%%%%%
\par\vspace{3mm}
%%%%%%%%%%%%%%%%
%      
%%%%%%%%%%%%%%%%%%%%%%%%%%%%%%%%%%%
In the case of $Y_{6,6}^{*,*}$, choose one of the two blocks 
of type $T_{\B}$ which it is $B_*$-adjacent to and discharge 
the weight $2$ of $Y_{6,6}$ to the block of type $T_{\B}$. 
In the case of $Z_{8,6}^{*,*}$ (resp. $Z_{10,6}^{bt,*}$), 
discharge $2$ out of its weight to each of the block of type 
$T_{\B}$ which its $f_8$ (resp. $f_{10}$) is $B_*$-adjacent 
to. Now, in each of the above $3$ cases and the first $7$ of 
$23$ cases, we have the situation that a block with its weight 
$0$ is $B_*$-adjacent to $T_{\B}$. We call that such blocks 
are type I as before. In each of the last $4$ of $23$ cases,
take a look at the face which is $B_*$-adjacent to only one 
$T_\B$ and discharge $2$ out of its weight $4$ to the block 
of type $T_{\B}$. In the case of 
$Y_{8,8}^{bt,bt}$, choose one of two faces of degree $8$ and 
discharge $2$ out of its weight $6$ to each of the blocks of 
type $T_{\B}$ which the face is $B_b$- or $B_t$-adjacent to. 
In the case of $Z_{10,8}^{bt,bt}$, discharge $2$ out of its 
weight $6$ to each of the blocks of type $T_{\B}$ which its 
$f_{10}$ is $B_b$- or $B_t$-adjacent to. Now, in each of the 
above $4$ cases and the second $7$ of $23$ cases, 
we have the situation that a block with its weight $2$ is 
$B_*$-adjacent to two blocks of type $T_{\B}$ with each weight 
$-2$ (if the two blocks are the same, we can discharge the 
weight $2$ to the block of type $T_{\B}$ and make its weight 
non-negative. Therefore, we do not think about such a case). 
We say that such blocks are type II as before.
%%%%%%%%%%%%%%%%%%%%%%%%%%%%%%%%%%%
%
%%%%%%%%%%%%%%%%
\par\vspace{3mm}
%%%%%%%%%%%%%%%%
%
%%%%%%%%%%%%%%%%%%%%%%%%%%%%%%%%%%%
Then we can constract paths as we did in Case 1-$\A$. Note that
Claim \ref{ClmArcs} holds in this case as well. Moreover, for 
each block of type I, we can find an arc $\x$ or $\y$ (resp.
$\y'$) in the boundary curve of $f_6$ which is $B_b$- (resp.
$B_t$-) adjacent to $T_{\B}$ from the primeness. Thus a 
similar claim to Claim \ref{Clm1aTpBtm} holds, which tells us 
that there does not exist such a path.
%%%%%%%%%%%%%%%%%%%%%%%%%%%%%%%%%%%
%
%%%%%%%%%%%%%%%%
\par\vspace{3mm}
%%%%%%%%%%%%%%%% 
%      
%%%%%%%%%%%%%%%%%%%%%%%%%%%%%%%%%%%%%%%%%%%%%%%%%% Case 1-C %%%%%
%
%%%%%%%%%%%%%%%%%%%%%%%%
\noindent{\bf Case 1-$\C$.}
%%%%%%%%%%%%%%%%%%%%%%%%
%
%%%%%%%%%%%%%%%%%%%%%%%%%%%%%%%%%%%
We show only the case that we have a block of type $\C$ as 
the previous case. Define a block of type $T_{\C}$ as before 
and take a look at the face $f$ which is $B_t$-adjacent to 
a block of type $T_{\C}$. Similarly to Case 1-$\B$, we have 
that $f$ cannot be $A$-adjacent to any block of type $T$ and 
that $f$ must be adjacent to at least one $f_2$. 
Therefore, we need to take a look at the blocks of type 
$Y_{i,j}^{p,q}$ or $Z_{i,j}^{p,q}$, where $p$, $q$ $\in$ 
$\{\cdot$, $t$, $a\}$ and one of them is $t$. It is easy to 
check that graph $G$ does not have $Y_6^t$, $Z_4^{\cdot}$, or 
$Z_6^a$. Additionally using Claim \ref{ClmAadj} and Claim 
\ref{Clm1Small}, we can see that the sum of the weights of 
such a block and the block of type $T_{\C}$ which the block is 
$B_t$-adjacent to is non-negative for each case. Therefore, we 
can discharge the weights of such blocks to blocks of type 
$T_{\C}$ to make the weight of every block non-negative.
%%%%%%%%%%%%%%%%%%%%%%%%%%%%%%%%%%%
%
%%%%%%%%%%%%%%%%
\par\vspace{3mm}                                                 
%%%%%%%%%%%%%%%%                   
%    
%%%%%%%%%%%%%%%%%%%%%%%%%%%%%%%%%%%%%%%%%%%%%%%%%% Case 1-G %%%%%
%
%%%%%%%%%%%%%%%%%%%%%%%%
\noindent{\bf Case 1-$\G$.}
%%%%%%%%%%%%%%%%%%%%%%%%
%
%%%%%%%%%%%%%%%%%%%%%%%%%%%%%%%%%%%
We show only the case that we have a block of type $\G$ as the 
previous case. Discharge the weight $2$ of the block to the
inner block of type $T$ which $\G$ is $A$-adjacent to (the one
whose boundary curve is surrounded by the boundary curve of the 
other on the diagram). Take a block of type $\G$ and the outer 
block of type $T$ and put them together. We call it a block of 
type $T_{\G}$ or simply $T_{\G}$, and so $w(T_{\G})$ $=$ $-2$. 
Take a look at the face $f$ which is $B_t$-adjacent to a block 
of type $T_{\G}$. Then, it cannot be adjacent to
any face of degree $2$, since its boundary curve is inside of 
the boundary curve of the face which is adjacent to the inner 
block of type $T$ at $v_d$. Therefore, we need to take a look 
at the blocks of type $X_i^p$, where $p$ is $t$ or $at$. It is 
easy to check that graph $G$ does not have $X_6^{at}$.
Additionally using Claim \ref{Clm1Small}, we can see that the 
sum of the weights of such a block and the block of type 
$T_{\G}$ which the block is $B_t$-adjacent to is non-negative 
for each case. Therefore, we can discharge the weights of such 
blocks to blocks of type $T_{\G}$ to make the weight of every 
block non-negative.
%%%%%%%%%%%%%%%%%%%%%%%%%%%%%%%%%%%
%
%%%%%%%%%%%%%%%%
\par\vspace{3mm}                                                 
%%%%%%%%%%%%%%%%                   
%      
%%%%%%%%%%%%%%%%%%%%%%%%%%%%%%%%%%%%%%%%%%%%%%%%%% Case 1-H %%%%%
%
%%%%%%%%%%%%%%%%%%%%%%%%
\noindent{\bf Case 1-$\H$.}
%%%%%%%%%%%%%%%%%%%%%%%%
%
%%%%%%%%%%%%%%%%%%%%%%%%%%%%%%%%%%%
We assume that we have a block of type $\H$. Define $T_{\H}$ as 
we did in Case 1-$\G$. Take a look at a face $f$ which is 
$B_t$-adjacent to $T_{\H}$. If $f$ is a face of a block $Y_{i,j}$ 
or $Z_{i,j}$, then let $f'$ be the other face of the block which 
is adjacent to a face of degree $2$ with $f$. Note that $f'$ 
cannot be $B_*$-adjacent to any blocks of type $T_{\H}$, since 
its boundary curve is inside of the boundary curve of the face 
of degee $8$ of $\H$. It is easy to see that $f$ is not 
a face of type $X_4^t$ or $X_6^{at}$ and that $f'$ is not a face
of degree $4$ or a face of type $Y_6^a$. Moreover, we can easily
obtain that graph $G$ does not have $Z_6^t$ or $Z_8^{at}$. 
Therefore, we can see that the sum of the weights of such a block 
and the block of type $T_{\H}$ which the block is $B_t$-adjacent 
to is non-negative for blocks of type $X_i^p$, $Y_{i,j}^{q,r}$, 
and $Z_{i,j}^{q,r}$, where $p$ and one of $q$, $r$ are of $\{t$, 
$at\}$ and the other is of $\{\cdot$, $a\}$. Therefore, we can 
discharge the weights of such blocks to blocks of type $T_{\H}$ 
to make the weight of every block non-negative.
%%%%%%%%%%%%%%%%%%%%%%%%%%%%%%%%%%%
%
%%%%%%%%%%%%%%%%
\par\vspace{3mm}                                                 
%%%%%%%%%%%%%%%%                   
%      
%%%%%%%%%%%%%%%%%%%%%%%%%%%%%%%%%%%%%%%%%%%%%%%%%% Case 1-N %%%%%
%
%%%%%%%%%%%%%%%%%%%%%%%%
\noindent{\bf Case 1-$N$.}
%%%%%%%%%%%%%%%%%%%%%%%%
%
%%%%%%%%%%%%%%%%%%%%%%%%%%%%%%%%%%%
Assume that we have at least one block of a type of $N$. We can 
similarly show the case that we have a block of a type of $N'$. 
We show this case step by step.
%%%%%%%%%%%%%%%%%%%%%%%%%%%%%%%%%%%
%
%%%%%%%%%%%%%%%%
\par\vspace{3mm} 
%%%%%%%%%%%%%%%%                   
%  
%%%%%%%%%%%%%%%%%%%%%%%%%%%%%%%%%%%
\noindent{\bf Step 1:}
To each of $\E$, $\F$, $\I$, and $\J$, discharge its weight $2$ 
to the inner block of type $T$ which it is $A$-adjacent to. 
Define $T_{\D}$ as the 
union of $T$ and $\D$. Define $T_{\E}$, $T_{\F}$, $T_{\I}$, and 
$T_{\J}$ as the union of the outer $T$ and $\E$, $\F$, $\I$, and 
$\J$, respectively. Call the block which is $B_b$-adjacent to 
$T_k$ a block of type $f_k$, where $k$ is $\D$, $\E$, $\F$, $\I$, 
or $\J$. No face can be $B_b$-adjacent to more than one such a 
block of type $T$ from Claim \ref{ClmCoexst} (ii). In addition, 
such a face can be adjacent to at most one face of degree $2$. 
Therefore, we may assume that the type of a block which has
such faces is $X_i^p$ or $Y_{i,j}^{q,r}$, where one of $q$, $r$ 
is $\cdot$ or $a$ and the other and $p$ are $b$ or $ab$. If 
the boundary curve of a block of type $T$ forms a diagram shown 
in Figure \ref{fig:GaDe} with boundary curves of some faces which 
it is adjacent to, then we say that the block of type 
$T$ is also type $\Ga$ or type $\De$. Then, we have the following.
%%%%%%%%%%%%%%%%%%%%%%%%%%%%%%%%%%%
%
%%%%%%%%%%%%%%%%%%%%%%%%%%%%%%%%%%%%%%%%%%% CLAIM Clm1nBadj %%%%%
\begin{Claim}\label{Clm1nBadj}\begin{itemize}
\item[(i)]   If $w(f_{\D})$ $+$ $w(T_{\D})$ $<$ $0$, then 
             $f_{\D}$ is type $X_6^{ab}$ or $Y_{8,4}^{ab,\cdot}$ 
             and the block of type $T$ which it is $A$-adjacent 
             to is type $\Ga$ or $\De$.
\item[(ii)]  If $w(f_{\E})$ $+$ $w(T_{\E})$ $<$ $0$, then $f_{\E}$ 
             is type $X_6^{ab}$ and the block of type $T$ which 
             it is $A$-adjacent to is type $\Ga$ or $\De$.
\item[(iii)] If $w(f_{\F})$ $+$ $w(T_{\F})$ $<$ $0$, then $f_{\F}$ 
             is type $X_6^{ab}$ and the block of type $T$ which 
             it is $A$-adjacent to is type $\Ga$.
\item[(iv)]  If $w(f_{\I})$ $+$ $w(T_{\I})$ $<$ $0$, then $f_{\I}$ 
             is type $X_6^{ab}$ and the block of type $T$ which 
             it is $A$-adjacent to is type $\Ga$.
\item[(v)]   $w(f_{\J})$ $+$ $w(T_{\J})$ $\geq$ $0$.
\end{itemize}\end{Claim}
%%%%%%%%%%%%%%%%%%%%%%%%%%%%%%%%%%%%%%%%%%%%%%%%% END CLAIM %%%%%
%
%%%%%%%%%%%%%%%%%%%%%%%%%%%%%%%%%%%%%%%%%%%%%%%%%%%%% PROOF %%%%%
\begin{proof}
(i) It is easy to see that $f_\D$ is $X_i^b$, $X_i^{ab}$, 
$Y_{i,4}^{b,\cdot}$, or $Y_{i,4}^{ab,\cdot}$. If $X_i^b$ $=$ 
$X_{\ge6}$, $X_i^{ab}$ $=$ $X_{\ge8}$, $Y_{i,4}^{b,\cdot}$ $=$ 
$Y_{\ge8,4}$, $Y_{i,4}^{ab,\cdot}$ $=$ $Y_{\ge10,4}$, then 
$w(f_\D)$ $+$ $w(T_\D)$ $\ge$ $0$. From Claim \ref{Clm1Small}, we 
do not have $X_4^b$. Let $v_d$$v_1$$v_2$$v_3$$v_4$$v_5$ be the 
boundary cycle of $Y_6^b$. Then, the face is adjacent to $f_2$ at 
$v_d$$v_1$. From the almost alternating property, it must be 
$B_b$-adjacent to $T_\D$ at $v_3$$v_4$. Then, we have that 
$b_5$ $=$ $b_\d$ or $b_\e$. The former case contradicts the 
reducedness and the latter case contradicts the minimality of 
$\t{L}$ (diagram V). Therefore, if $w(f_{\D})$ $+$ 
$w(T_{\D})$ $<$ $0$, then $f_{\D}$ is $X_6^{ab}$ or 
$Y_{8,4}^{ab,\cdot}$. In the first case, let 
$v_1$$v_2$$v_3$$v_4$$v_5$$v_6$ be the boundary cycle of 
$X_6^{ab}$ and let $X_6^{ab}$ be $A$-adjacent to $T'$ at $b_1$ 
$=$ $b_{\c'}$. From the almost alternating property, we have that 
$b_2$ $=$ $b_\a$ or $b_4$ $=$ $b_\a$. Since the length of curve 
$b_1$$b_2$$b_3$$b_4$$b_5$$b_6$ is $6$, we have that $b_2$ $=$ 
$b_\a$ and $b_4$ $=$ $b_d$. Then we have that $b_5$ $=$ $b_\d$,
$b_\e$ or that $b_6$ $=$ $b_\d$, $b_\e$. The first and the fourth
cases contradict the reducedness. In the second and the third
case, we obtain the diagram of type $\Ga$ and type $\De$, 
respectively. For the case of $Y_{8,4}^{ab,\cdot}$,
let $v_d$$v_1$$v_2$$v_3$$v_4$$v_5$$v_6$$v_7$ be the 
boundary cycle of $Y_{8,4}$ and let $Y_{8,4}$ be adjacent to 
$f_2$ at $v_d$$v_1$. Then, we do not have that $b_\b$ $=$ $b_3$ 
or $b_5$, since the length of curve 
$b_d$$b_1$$b_2$$b_3$$b_4$$b_5$$b_6$$b_7$ is $8$ and it must pass 
arc $b_\b$$b_\a$ and bubble $b_\c$. If $b_3$ $=$ $b_\b$, then we 
have that $b_5$ $=$ $b_{\c'}$ considering the length of the curve. 
In the former case, we have that $b_6$ $=$ $b_\d$, $b_6$ $=$ $b_\e$, 
$b_7$ $=$ $b_\d$, or $b_7$ $=$ $b_\e$. The second and the third 
cases contradict the mimimality of $\t{L}$ from Lemma 
\ref{LemMinimality}. In the first and the 
fourth case, we obtain $T$ of type $\De$ and $\Ga$, respectively.
(ii) $\sim$ (iv) Note that $T_\E$ is $T$ of type $\Ga$, and 
$T_\F$ and $T_\I$ are type $\De$. Now it is easy to see that those 
statements hold. (v) We can see this from the primeness and the 
reducedness of $\t{L}$. 
\end{proof}
%%%%%%%%%%%%%%%%%%%%%%%%%%%%%%%%%%%%%%%%%%%%%%%%% END PROOF %%%%%
%
%%%%%%%%%%%%%%%%%%%%%%%%%%%%%%%%%%%
\par\vspace{3mm}\begin{figure}[htbp]\caption{$\Ga$ and $\De$}
\label{fig:GaDe}\end{figure}\par\vspace{3mm}
%%%%%%%%%%%%%%%%%%%%%%%%%%%%%%%%%%%
% 
%%%%%%%%%%%%%%%%%%%%%%%%%%%%%%%%%%%
If $w(f_k)$ $+$ $w(T_k)$ $\ge$ $0$, then discharge $2$ out of
the weight of $w(f_k)$ to the block of type $T_k$, where $k$ is 
$\D$, $\E$, $\F$, $\I$, or $\J$. If every such a sum is 
non-negative, then we are done. If there is a case that the sum 
is negative (we know that it is $-2$ from the above claim), then 
we go on to the next step.
%%%%%%%%%%%%%%%%%%%%%%%%%%%%%%%%%%%
%
%%%%%%%%%%%%%%%%
\par\vspace{3mm} 
%%%%%%%%%%%%%%%%                   
%  
%%%%%%%%%%%%%%%%%%%%%%%%%%%%%%%%%%%
\noindent{\bf Step 2:}
Take weight $2$ of the block of type $T$ which the block of type 
$f_k$ is $A$-adjacent to (thus its weight goes dowm from $0$ to
$-2$) and give it to the block of type $T_k$, where $k$ is $\D$, 
$\E$, $\F$, or $\I$. We can see from the proof of the previous 
claim that the boundary curves of the bottom faces of all the 
blocks of type $T$ which the blocks of type $f_k$ are $A$-adjacent 
to pass the same four bubbles and that such a block of type $T$
is also type $\Ga$ or $\De$. Now take a look at a face which
is $B_b$-adjacent to a block of type $\Ga$ or $\De$. Following 
the proof of Claim \ref{Clm1nBadj} (ii) and (iii), the 
only case that the sum of the weights of the block and blocks
of type $\Ga$ or $\De$ which it is $B_b$-adjacent to is negative 
is that the block is $X_6^{ab}$ and then we obtain a block of 
type $\Ga$ or $\De$ again, which the block is $A$-adjacent to.
%%%%%%%%%%%%%%%%%%%%%%%%%%%%%%%%%%%
%
%%%%%%%%%%%%%%%%
\par\vspace{3mm} 
%%%%%%%%%%%%%%%%                   
%  
%%%%%%%%%%%%%%%%%%%%%%%%%%%%%%%%%%%
\noindent{\bf Step 3:}
Now, we are at the beginning of Step 2. Since $G$ is a finite 
graph, we can finally reach the situation so that we can 
discharge weight of a block to a block of type $\Ga$ or $\De$ 
with negative weight by continuing this process.
%%%%%%%%%%%%%%%%%%%%%%%%%%%%%%%%%%%
%
%%%%%%%%%%%%%%%%
\par\vspace{5mm}                                                 
%%%%%%%%%%%%%%%%    
%
%

%%%%%%%%%%%%%%%%%
%     
%%%%%%%%%%%%%%%%%%%%%%%%%%%%%%%%%%%%%%%%%%%%%%%% Case 2 & 3%%%%%
%
%%%%%%%%%%%%%%%%%
%
%%%%%%%%%%%%%%%%%%%%%%%%%%%%%%%%%%%     
\begin{center}{\bf Case 2.}\end{center}
%%%%%%%%%%%%%%%%%%%%%%%%%%%%%%%%%%%     
%
%%%%%%%%%%%%%%%%
\par\vspace{5mm}
%%%%%%%%%%%%%%%%
%        
%%%%%%%%%%%%%%%%%%%%%%%%%%%%%%%%%%%     
In this case, we look at faces which are $C_*$-adjacent to 
blocks of type $U$, where we use the notation $C_*$ to mean 
$C_l$ or $C_r$. From Claim \ref{ClmBdTU}, we can see that 
every face can be $C_*$-adjacent to at most one block of type
$U$ at most once. Then we have $7$ types of blocks which are
$C_*$-adjacent to blocks of type $U$; $X_i^*$, $Y_{i,j}^{p,q}$, 
and $Z_{i,j}^{p,q}$ with $\{p,q\}$ $=$ $\{\cdot,*\}$, 
$\{*,\cdot\}$, or $\{*,*\}$, where $*$ stands for $l$ or $r$. 
Here, we have the following claim.
%%%%%%%%%%%%%%%%%%%%%%%%%%%%%%%%%%%
%
%%%%%%%%%%%%%%%%
\par\vspace{3mm}              
%%%%%%%%%%%%%%%%                   
%     
%%%%%%%%%%%%%%%%%%%%%%%%%%%%%%%%%%%%%%%%%%%%% CLAIM ClmCadj %%%%%
\begin{Claim}\label{ClmCadj}\begin{itemize}
\item[(i)]   No face of degree $4$ can be $C_*$-adjacent to a
             block of type $U$.
\item[(ii)]  No face of degree $4$ can be adjacent to two 
             faces of degree $2$ with any face of degree $6$
             which is $C_*$-adjacent to a block of type $U$.     
\item[(iii)] No face of degree $6$ can be adjacent to two 
             faces of degree $2$ with any other face of degree 
             $6$ which is $C_*$-adjacent to a block of type $U$. 
\end{itemize}\end{Claim}
%%%%%%%%%%%%%%%%%%%%%%%%%%%%%%%%%%%%%%%%%%%%%%%%% END CLAIM %%%%%
%
%%%%%%%%%%%%%%%%%%%%%%%%%%%%%%%%%%%%%%%%%%%%%%%%%%%%% PROOF %%%%%
\begin{proof}
(i) Let $v_\f$$v_\g$$v_1$$v_2$ be the boundary cycle of the face
of degree $4$. If $v_1$ $=$ $v_d$, then it contradicts the 
minimality of $\t{L}$ (diagram VI). If $v_2$ 
$=$ $v_d$, then there exists a face which has two black vertices
on its boundary, which contradicts Lemma \ref{LemEveryCurve}.
%%%%%%%%%%%%%%%%%%%%%%%%%%%%%%%%%%%        
(ii), (iii) Let $v_\f$$v_\g$$v_1$$v_2$$v_3$$v_4$ be the 
boundary cycle of the face of degree $6$. Similarly to the
previous case, we can show the cases of that $v_1$ $=$ $v_d$ 
and $v_4$ $=$ $v_d$. If $v_3$ $=$ $v_d$, then it also 
contradicts the minimality of $\t{L}$ (diagram VII). 
Therefore, we have that $v_2$ $=$ $v_d$.
Now let $v_1$$v_2$$v_3$$v_5$ and $v_1$$v_2$$v_3$$v_6$$v_7$$v_8$
be the boundary curve of the face of degree $4$ and $6$ in the
statement, respectively. Then, note that each bounday curve is
surrounded by that of the face of degree $6$ which is 
$C_*$-adjacent to a block of type $U$. If $b_\f$ $=$ $b_5$, 
$b_6$, or $b_8$, then it contradicts the minimality of $\t{L}$ 
(diagram VI or VII). If $b_7$ $=$ $b_\f$, then it
contradicts the primeness.
\end{proof}
%%%%%%%%%%%%%%%%%%%%%%%%%%%%%%%%%%%%%%%%%%%%%%%%% END PROOF %%%%%
%
%%%%%%%%%%%%%%%%
\par\vspace{3mm}
%%%%%%%%%%%%%%%%
%      
%%%%%%%%%%%%%%%%%%%%%%%%%%%%%%%%%%%        
If $X_i$ is $C_l$- and $C_r$-adjacent to a block of type $U$, 
then we call it $X_i^l$ and $X_i^r$, respectively. Since we do 
not have $X_4^*$ from Claim \ref{ClmCadj}, $w(X_i^*)=$ $i-4\ge$ 
$6-4=$ $2$. Similarly we have that $Y_{i,j}^{*,\cdot}$ $=$ 
$Y_{\ge6,\ge4}$, $Y_{i,j}^{\cdot,*}$ $=$ $Y_{\ge4,\ge6}$, 
$Y_{i,j}^{*,*}$ $=$ $Y_{\ge6,\ge6}$, $Z_{i,j}^{*,\cdot}$ $=$
$Z_{\ge8,\ge4}$, $Z_{i,j}^{\cdot,*}$ $=$ $Z_{\ge8,\ge6}$, and
$Z_{i,j}^{*,*}$ $=$ $Z_{\ge8,\ge6}$. Therefore we have that
$w(Y_{i,j}^{*,\cdot})$ $=$ $i$ $+$ $j$ $-$ $10$ $\ge$ $0$,
$w(Y_{i,j}^{\cdot,*})$ $\ge$ $0$, $w(Y_{i,j}^{*,*})$ $\ge$ $2$,
$w(Z_{i,j}^{*,\cdot})$ $=$ $i$ $+$ $j$ $-$ $12$ $\ge$ $0$,
$w(Z_{i,j}^{\cdot,*})$ $\ge$ $2$, and $w(Z_{i,j}^{*,*})$ $\ge$ 
$2$. Then, for each block which is $C_*$-adjacent to blocks of 
type $U$, discharge $2$ out of its weight to each of the blocks 
of type $U$ if the sum of the weights of the block and all the 
blocks of type $U$ is non-negative. The type of block such 
that the sum is negative is $Y_{6,4}^{*.\cdot}$, 
$Y_{4,6}^{\cdot,*}$, $Y_{6,6}^{*,*}$, $Z_{8,4}^{*,\cdot}$, or 
$Z_{8,6}^{*,*}$. For each of these blocks, we say that it is 
type II if it is $C_l$- and $C_r$-adjacent to two blocks of 
type $U$. Otherwise, we say that it is type I. Then, we have 
the following claim, where we say that a block of type $U$ is 
$D_l$- or $D_r$-adjacent to a face if the boundary curves of 
the block of type $U$ and the block containing the face form 
diagram $\Th$ shown in Figure \ref{fig:Th}.
%%%%%%%%%%%%%%%%%%%%%%%%%%%%%%%%%%%
%
%%%%%%%%%%%%%%%%%%%%%%%%%%%%%%%%%%%
\begin{figure}[htbp]\caption{$\Th$}
\label{fig:Th}\end{figure}\par\vspace{3mm}
%%%%%%%%%%%%%%%%%%%%%%%%%%%%%%%%%%%                 
%     
%%%%%%%%%%%%%%%%%%%%%%%%%%%%%%%%%%%%%%%%%%% CLAIM Clm2TypeI %%%%%
\begin{Claim}\label{Clm2TypeI} For every block of type I, the 
block of type $U$ which it is $C_l$- $($resp. $C_r$-$)$ adjacent 
to is $D_l$- $($resp. $D_r$-$)$ adjacent to a face of the block.  
\end{Claim}
%%%%%%%%%%%%%%%%%%%%%%%%%%%%%%%%%%%%%%%%%%%%%%%%% END CLAIM %%%%%
%
%%%%%%%%%%%%%%%%%%%%%%%%%%%%%%%%%%%%%%%%%%%%%%%%%%%%% PROOF %%%%%
\begin{proof}
The block of type I is $Y_{6,4}^{*,\cdot}$, $Y_{4,6}^{\cdot,*}$,
$Y_{6,6}^{l,l}$, $Y_{6,6}^{r,r}$, $Z_{8,4}^{*,\cdot}$,
$Z_{8,6}^{l,l}$, or $Z_{8,6}^{r,r}$. We show only the cases of
$Y_{6,4}^{l,\cdot}$, $Y_{6,6}^{l,l}$, $Z_{8,6}^{l,l}$.
%%%%%%%%%%%%%%%%%%%%%%%%%%%%%%%%%%%
%
%%%%%%%%%%%%%%%%%%%%%%%%%%%%%%%%%%%
In the case of $Y_{6,4}^{l,\cdot}$, let 
$v_\f$$v_\g$$v_1$$v_2$$v_3$$v_4$ be the boundary curve of 
$Y_6^l$. We can show the cases that $v_1$ $=$ $v_d$ 
and $v_4$ $=$ $v_d$ similarly to Claim \ref{ClmCadj} (i).
Assume that $v_2$ $=$ $v_d$. Considering the face of degree
$4$ of $Y_{6,4}$, we can see that it contradicts the 
minimality (VI or VII) or primeness of $\t{L}$ in both cases 
that two faces of $Y_{6,4}$ are adjacent to a face of degree
$2$ at $v_1$$v_2$ and $v_2$$v_3$. Next assume that $v_3$ 
$=$ $v_d$. Then the two faces of $Y_{6,4}$ are adjacent to 
a face of degree $2$ at $v_2$$v_3$ or $v_3$$v_4$. The latter
case contradicts the minimality of $\t{L}$ (diagram VII). In 
the former case, let $v_2$$v_3$$v_5$$v_6$ be the boundary cycle 
of the face of degree $4$. Then, we have that $b_5$ $=$ $b_\g$, 
$b_5$ $=$ $b_\h$, $b_6$ $=$ $b_\g$, or $b_6$ $=$ $b_\h$. The
first two cases contradict the primeness and the last case
contradicts the minimality of $\t{L}$ (diagram VI). In the third
case, we obtain diagram $\Th$.
%%%%%%%%%%%%%%%%%%%%%%%%%%%%%%%%%%%
%
%%%%%%%%%%%%%%%%
\par\vspace{3mm}      
%%%%%%%%%%%%%%%%                   
%     
%%%%%%%%%%%%%%%%%%%%%%%%%%%%%%%%%%%
In the case of $Y_{6,6}^{l,l}$, considering the face of degree
$6$ whose boundary curve passes the leftside of the dealternator
on the diagram and following the previous case, we have that
$v_3$ $=$ $v_d$ and the two faces of degree $6$ are adjacent to
a face of degree $2$ at $v_2$$v_3$. Now let 
$v_2$$v_3$$v_5$$v_6$$v_7$$v_8$ be the boundary curve of the
other face of degree $6$. Since it is also $C_l$-adjacent to 
a block of type $U$, we have that $b_6$ or $b_8$ $=$ $b_\f$ 
from the almost alternating property. The former case 
contradicts the primeness. In the latter case, we obtain 
diagram $\Th$.
%%%%%%%%%%%%%%%%%%%%%%%%%%%%%%%%%%%
%
%%%%%%%%%%%%%%%%
\par\vspace{3mm}            
%%%%%%%%%%%%%%%%                   
%     
%%%%%%%%%%%%%%%%%%%%%%%%%%%%%%%%%%%
In the case of $Z_{8,6}^{l,l}$, considering the face of degree
$6$ and following the proof of Claim \ref{ClmCadj} (ii) and 
(iii), we have that $v_2$ $=$ $v_d$. Now let 
$v_1$$v_2$$v_3$$v_5$$v_6$$v_7$$v_8$$v_9$ be the boundary cycle
of the face of degree $8$. Since the face is also $C_l$-adjacent 
to a block of type $U$, we have that $b_6$ or $b_8$ $=$ $b_\g$ 
considering the almost alternating property. The former case 
contradicts the mimimality of $\t{L}$ (diagram VII). In the 
latter case, we obtain diagram $\Th$.
\end{proof}
%%%%%%%%%%%%%%%%%%%%%%%%%%%%%%%%%%%%%%%%%%%%%%%%% END PROOF %%%%%
% 
%%%%%%%%%%%%%%%%%%%%%%%%%%%%%%%%%%%
For each of $Y_{6,6}^{*,*}$ and $Z_{6,8}^{*,*}$ of type I, 
discharge its weight $2$ to the block of type $U$ which is 
not $D_l$- or $D_r$-adjacent to any face of the block. Then, 
we may conclude that if we still have a block of type $U$ 
with negative weight, then it is $D_l$- (resp. $D_r$-) 
adjacent to a block with weight $0$, or it is $C_l$- 
(resp. $C_r$-) adjacent to a block with weight $2$ which is 
$C_r$- (resp. $C_l$-) adjacent to another block of type $U$ 
with negative weight. Then, we can constract finite paths by 
regarding blocks of type $U$ and blocks of type I and II as 
edges and vertices, respectively. However then, clearly from 
their diagrams, if there exists a block of type $U$ which is 
$D_l$-adjacent to a face of a block of type I, then there does 
not exist any block of type $U$ which is $D_r$-adjacent to a 
face of a block of type I. Therefore, there does not exist 
such a path, since its ends should come from blocks of type I.
%%%%%%%%%%%%%%%%%%%%%%%%%%%%%%%%%%%
%
%%%%%%%%%%%%%%%%
\par\vspace{5mm}              
%%%%%%%%%%%%%%%%  
% 
%
%%%%%%%%%%%%%%%%%%%%%%%%%%%%%%%%%%%     
\begin{center}{\bf Case 3.}\end{center}
%%%%%%%%%%%%%%%%%%%%%%%%%%%%%%%%%%%     
%
%%%%%%%%%%%%%%%%
\par\vspace{5mm}
%%%%%%%%%%%%%%%%
%        
%%%%%%%%%%%%%%%%%%%%%%%%%%%%%%%%%%%    
Now we have a block of type $T$ and a block of type $U$. If 
$b_\e$ $=$ $b_\g$ or $b_\a$ $=$ $b_\h$, then it contradicts 
the minimality of $\t{L}$ (diagram III or VI). If $b_\a$ $=$ 
$b_\g$, then it contradicts the primeness. And if $b_\e$ $=$ 
$b_\h$ and $b_\d$ $\ne$ $b_\g$ (or $b_\e$ $\ne$ $b_\h$ and 
$b_\d$ $=$ $b_\g$), then it also contradicts the primeness. 
If $b_\e$ $=$ $b_\h$ and $b_\d$ $=$ $b_\g$, then we treat the 
block of type $T$ as a block of type $U$. In the case that 
$v_2$ $=$ $v_d$ of Claim \ref{ClmCadj} (i), we obtain a 
non-reduced diagram. If we have that $b_\e$ $=$ $b_\h$ and 
$b_\d$ $=$ $b_\g$ for every block of type $T$, then we are 
done by following Case 2. Therefore, we may assume that there 
exists at least one block of type $T$ such that $b_\d$ $\ne$ 
$b_\g$ and $b_\e$ $\ne$ $b_\h$. Since the diagram shown in 
Figure \ref{fig:TwithU} and its mirror image cannot coexist, 
we may assume that $T$ and $U$ replaced on the diagram are 
shown in Figure \ref{fig:TwithU}.
%%%%%%%%%%%%%%%%%%%%%%%%%%%%%%%%%%%
%
%%%%%%%%%%%%%%%%%%%%%%%%%%%%%%%%%%%
First of all, take a look at the blocks which are $A$-adjacent to
blocks of type $T$ and discharge $2$ out of its weight to each of 
the blocks of type $T$ if the sum of the weights of the block and 
all the blocks of type $T$ is non-negative. After discharging, we
have the two cases; there are no blocks of type $T$ with negative 
weight ({\bf Case 3-1}) and there exists a block of type $T$ with 
negative weight ({\bf Case 3-2}). 
%%%%%%%%%%%%%%%%%%%%%%%%%%%%%%%%%%%
%
%%%%%%%%%%%%%%%%%%%%%%%%%%%%%%%%%%%
\par\vspace{3mm}\begin{figure}[htbp]\caption{$T$ with $U$}
\label{fig:TwithU}\end{figure}\par\vspace{3mm}
%%%%%%%%%%%%%%%%%%%%%%%%%%%%%%%%%%%
%
%%%%%%%%%%%%%%%%%%%%%%%%%%%%%%%%%%%
\begin{center}{\bf Case 3-1.}\end{center}\par\vspace{3mm}
%%%%%%%%%%%%%%%%%%%%%%%%%%%%%%%%%%%
%
%%%%%%%%%%%%%%%%%%%%%%%%%%%%%%%%%%%
Take a look at the blocks which are $C_*$-adjacent to blocks of 
type $U$. Now we have the following claim.
\par\vspace{3mm}
%%%%%%%%%%%%%%%%%%%%%%%%%%%%%%%%%%%
%
%%%%%%%%%%%%%%%%%%%%%%%%%%%%%%%%%%%%%%%%%%% CLAIM Clm31Small %%%%%
\begin{Claim}\label{Clm31Small} Graph $G$ does not have any
blocks of type $X_4^*$, $Y_4^{\cdot}$, $Z_6^*$, $Z_{6,i}^{a,*}$. 
\end{Claim}
%%%%%%%%%%%%%%%%%%%%%%%%%%%%%%%%%%%%%%%%%%%%%%%%% END CLAIM %%%%%
%
%%%%%%%%%%%%%%%%%%%%%%%%%%%%%%%%%%%     
\begin{proof}
From Claim \ref{ClmCadj}, it is sufficient to show only the last
three cases. If there exists a block of type $Y_4^{\cdot}$, then
it contradicts the primeness, the reducedness, or the minimality 
of $\t{L}$ (diagram VI or VII). Let $v_d$$v_1$$v_2$$v_3$$v_4$$v_5$ 
be the boundary cycle of $Z_6^*$. Cosidering the length of the 
curve, we can assume that $b_1$ $=$ $b_\a$, $b_1$ $=$ $b_\b$, $b_5$
$=$ $b_\d$, or $b_5$ $=$ $b_\e$. Then, it contradicts the primeness 
or the minimality of $\t{L}$ (diagram VI or VII). Take a look at a 
block of type $Z_6^a$ and let $v_\b$$v_\c$$v_\d$$v_1$$v_2$$v_3$ be 
its boundary cycle. From the minimality of $\t{L}$ (diagram VI or 
VII), we obtain that $v_2$ $=$ $v_d$ and thus the curve surrounds 
the boundary curve of the face $f$ which is adjacent to two faces 
of degree $2$ with $Z_6^a$. Therefore, face $f$ cannot be 
$C_*$-adjacent to a block of type $U$.
\end{proof}
%%%%%%%%%%%%%%%%%%%%%%%%%%%%%%%%%%%     
%
%%%%%%%%%%%%%%%%
\par\vspace{3mm}              
%%%%%%%%%%%%%%%%                   
% 
%%%%%%%%%%%%%%%%%%%%%%%%%%%%%%%%%%%
Note that no face can be $A$-adjacent to $T$ and $C_*$-adjacent 
to $U$. Thus we have $11$ types of blocks which are
$C_*$-adjacent to blocks of type $U$; $X_i^*$, $Y_{i,j}^{p,q}$, 
and $Z_{i,j}^{p,q}$ with $\{p,q\}$ $=$ $\{\cdot,a,*\}$ and
$p$ or $q$ $=$ $*$. 
%%%%%
%From Claim \ref{Clm31Small}, we have that
%$X_i^*$ $=$ $X_{\ge6}$,
%$Y_{i,j}^{*,\cdot}$ $=$ $Y_{\ge6,\ge6}$, 
%$Y_{i,j}^{\cdot,*}$ $=$ $Y_{\ge6,\ge6}$, 
%$Y_{i,j}^{*,a}$ $=$ $Y_{\ge6,\ge6}$,
%$Y_{i,j}^{a,*}$ $=$ $Y_{\ge6,\ge6}$,
%$Y_{i,j}^{*,*}$ $=$ $Y_{\ge6,\ge6}$,
%$Z_{i,j}^{\cdot,*}$ $=$ $Z_{\ge8,\ge8}$,
%$Z_{i,j}^{*,\cdot}$ $=$ $Z_{\ge8,\ge6}$,
%$Z_{i,j}^{*,a}$ $=$ $Z_{\ge8,\ge8}$ or $Z_{\ge10,\ge6}$,
%$Z_{i,j}^{a,*}$ $=$ $Z_{\ge8,\ge8}$, and
%$Z_{i,j}^{*,*}$ $=$ $Z_{\ge8,\ge8}$.
%%%%%
Then, for each block which is $C_*$-adjacent to blocks of 
type $U$, discharge $2$ out of its weight to each of the blocks 
of type $U$ if the sum of the weights of the block and all the 
blocks of type $U$ is non-negative. From Claim \ref{Clm31Small},
The type of block such that the sum is negative is 
$Y_{6,6}^{*,a}$, $Y_{6,6}^{a,*}$, or $Y_{6,6}^{*,*}$. 
\par\vspace{3mm}
%%%%%%%%%%%%%%%%%%%%%%%%%%%%%%%%%%%
%
%%%%%%%%%%%%%%%%%%%%%%%%%%%%%%%%%%%
Now let us take a look at a block of type $Y_6^*$. First, let
$v_\f$$v_\g$$v_1$$v_2$$v_3$$v_4$ be the boundary cycle of 
$Y_6^l$. Then we have that $v_1$ $\ne$ $v_d$ and $v_4$ $\ne$ 
$v_d$ from the proof of Claim \ref{ClmCadj} (i). Assume that 
$v_2$ $=$ $v_d$. If $b_1$ $=$ $b_\d$, then it contradicts the 
minimality of $\t{L}$ (diagram VI). Thus we have that $b_1$ 
$=$ $b_\e$ and then the block is adjacent to a face of degree 
$2$ at $v_2$$v_3$, otherwise it contradicts the primeness. 
Therefore we have that $b_3$ $\ne$ $b_\a$, $b_3$ $\ne$ $b_\b$,
and $b_4$ $\ne$ $b_\a$ from the primeness and the minimality 
of $\t{L}$ (diagram VII). Thus we have that $b_4$ $=$ $b_\b$ 
and then we obtain diagram $D_1$ in Figure \ref{fig:y6landy6r}. 
Similary we obtain $D_2$ in the case that $v_3$ $=$ $v_d$. Next, 
let $v_\h$$v_\g$$v_1$$v_2$$v_3$$v_4$ be the boundary cycle of 
$Y_6^r$. Following the previous case, we see that the case that 
$v_2$ $=$ $v_d$ contradicts the primeness and we obtain diagram 
$D_3$ in the case that $v_3$ $=$ $v_d$.
%%%%%%%%%%%%%%%%%%%%%%%%%%%%%%%%%%%
%
%%%%%%%%%%%%%%%%%%%%%%%%%%%%%%%%%%%
\par\vspace{3mm}\begin{figure}[htbp]
\caption{$Y_6^l$ with $Y_6^r$}
\label{fig:y6landy6r}\end{figure}\par\vspace{3mm}
%%%%%%%%%%%%%%%%%%%%%%%%%%%%%%%%%%%
% 
%%%%%%%%%%%%%%%%%%%%%%%%%%%%%%%%%%%
Then we can see that we do not have a block of type $Y_{6,6}^{l,l}$
from diagrams $D_1$ and $D_2$ paying attention to the boundary 
curves of faces of degree $2$. And we also do not have a block of
type $Y_{6,6}^{r,r}$, since the boundary curve of $Y_6^r$ must
pass the rightside of the dealternator on the diagram. In addition,
note that the diagrams of $Y_6^l$ of $Y_{6,6}^{a,*}$ and 
$Y_{6,6}^{*,a}$ must be $D_1$, since $b_\c$ is surrounded by the
boundary curve of $Y_6^l$ in diagram $D_2$. Now we can constract 
finite paths by regarding the blocks of type $U$ as edges and 
$Y_{6,6}^{r,l}$, $Y_{6,6}^{l,r}$, $Y_{6,6}^{a,*}$, and 
$Y_{6,6}^{*,a}$ as vertices. Then, their ends must come from 
$Y_{6,6}^{*,a}$ or $Y_{6,6}^{a,*}$ and thus we have diagrams $D_1$
and $D_3$. However, we can see that it does not happen from Figure
\ref{fig:y6landy6r}. Therefore we can conclude that we do not have
such paths.\par\vspace{3mm}
%%%%%%%%%%%%%%%%%%%%%%%%%%%%%%%%%%%
%
%%%%%%%%%%%%%%%%%%%%%%%%%%%%%%%%%%%
\begin{center}{\bf Case 3-2.}\end{center}\par\vspace{3mm}
%%%%%%%%%%%%%%%%%%%%%%%%%%%%%%%%%%%
%
%%%%%%%%%%%%%%%%%%%%%%%%%%%%%%%%%%%
In this case, we see that the boundary curves of the bottom 
faces of blocks of type $T$ with negative weight pass the same 
four bubbles from the following claim and Claim \ref{ClmCoexst}. 
In the rest of this case, we use only this fact and we do not 
care about the type of a block which is $A$-adjacent to $T$.
%%%%%%%%%%%%%%%%%%%%%%%%%%%%%%%%%%%
%
%%%%%%%%%%%%%%%%%%%%%%%%%%%%%%%%%%%%%%% CLAIM Clm3Situation %%%%%
\begin{Claim}\label{Clm3Situation} After discharging, the weight 
of any block of type $T$ is non-negative, or $G$ contains $\A$, 
$\G^*$, or $\H$. \end{Claim}
%%%%%%%%%%%%%%%%%%%%%%%%%%%%%%%%%%%%%%%%%%%%%%%%% END CLAIM %%%%%
%
%%%%%%%%%%%%%%%%%%%%%%%%%%%%%%%%%%%%%%%%%%%%%%%%%%%%% PROOF %%%%%
\begin{proof}
The diagram obtained from any block of $N'$ and $N$ contains 
$S$ and $S'$, respectively. It contradicts the minimality of 
$\t{L}$ (diagram VI or VII), since we have $U$ as well. The 
diagram obtained from $\B$ (resp. $\C$) contains an arc 
connecting $b_\a$ and $b_\e$ (resp. $b_\b$ and $b_\e$). However, 
since $b_\a$ (resp. $b_\b$) is in one of the two regions of 
$S^2$ $-$ $b_d$$b_\f$$b_\g$$b_\h$ and $b_\e$ is in the other, 
no arc can connect them.
\end{proof}
%%%%%%%%%%%%%%%%%%%%%%%%%%%%%%%%%%%%%%%%%%%%%%%%% END PROOF %%%%%
%
%%%%%%%%%%%%%%%%
\par\vspace{3mm}              
%%%%%%%%%%%%%%%%                   
% 
%%%%%%%%%%%%%%%%%%%%%%%%%%%%%%%%%%%
Take a look at blocks which are $B_b$-adjacent to $T$ or 
$C_*$-adjacent to $U$. Then we have the following.
%%%%%%%%%%%%%%%%%%%%%%%%%%%%%%%%%%%
%
%%%%%%%%%%%%%%%%%%%%%%%%%%%%%%%%%%%%%%%%%% CLAIM Clm32Small %%%%%
\begin{Claim}\label{Clm32Small} Graph $G$ does not have any
block of type $X_4^b$, $X_6^b$, $X_8^{b*}$, $Z_8^{ab}$, 
$Y_{8,6}^{ab,*}$, $Y_{8,6}^{ab,a}$, $Y_{8,8}^{ab,ab}$, 
or $Z_{i,6}^{b,a}$.
\end{Claim}
%%%%%%%%%%%%%%%%%%%%%%%%%%%%%%%%%%%%%%%%%%%%%%%%% END CLAIM %%%%%
%
%%%%%%%%%%%%%%%%%%%%%%%%%%%%%%%%%%%%%%%%%%%%%%%%%%%%% PROOF %%%%%
\begin{proof}
We omit the proof for the first three and the last cases. Take 
the boundary cycle $v_\d$$v_\c$$v_\b$$v_1$$v_2$$v_3$$v_4$$v_5$ 
of $Y_8^{ab}$. From the almost alternating property and the 
length of the curve, $Y_8^{ab}$ is $B_b$-adjacent to $T$ at 
$v_\b$$v_1$ or $v_2$$v_3$. 
Since the first case contradicts the minimality of $\t{L}$
(diagram VII), we may consider the second case. Then we have that 
$v_5$ $=$ $v_d$ from the minimality of $\t{L}$ (diagram VII).
Then, we can see that the boundary curve passes the leftside
of the dealternator. Therefore, we do not have $Y_{8,8}^{ab,ab}$.
Considering the diagram, it is also easy to see that we do not
have $Z_8^{ab}$, $Y_{8,6}^{ab,*}$ or $Y_{8,6}^{ab,a}$.
\end{proof}
%%%%%%%%%%%%%%%%%%%%%%%%%%%%%%%%%%%%%%%%%%%%%%%%% END PROOF %%%%%
%
%%%%%%%%%%%%%%%%
\par\vspace{3mm}              
%%%%%%%%%%%%%%%%                   
% 
%%%%%%%%%%%%%%%%%%%%%%%%%%%%%%%%%%%
We know that no face can be $A$-adjacent to $T$ and $C_*$-adjacent 
to $U$, that any face can be $C_*$-adjacent to at most one $U$ at 
most once and, from the above fact, that any face can be 
$B_b$-adjacent to at most one $T$ at most once. Thus we have $68$ 
types of blocks which are $B_b$-adjacent to $T$ or $C_*$-adjacent 
to $U$; $X_i^p$, $Y_{i,j}^{q,r}$, and $Z_{i,j}^{q,r}$ with 
$\{q,r\}$ $=$ $\{\cdot,a,*,b,ab,*b\}$ and $p$ and one of $q$ and 
$r$ are of $\{*,b,*b,ab\}$. Then, for each block which is 
$C_*$-adjacent to a block of type $U$ or $B_b$-adjacent to a 
block of type $T$, discharge $2$ out of its weight to each 
of the blocks of type $T$ and $U$ if the sum of the weights of 
the block and all the blocks of type $T$ and $U$ is non-negative.
From Claim \ref{Clm31Small} and Claim \ref{Clm32Small}, we have 
that the type of a block such that the sum is negative is
$Y_{6,6}^{*,a}$, $Y_{6,6}^{a,*}$, or $Y_{6,6}^{*,*}$. However, 
then we can conclude that the sum of the weights of all faces is
non-negative following the previous case.
%%%%%%%%%%%%%%%%%%%%%%%%%%%%%%%%%%%
%
%%%%%%%%%%%%%%%%
\par\vspace{5mm}
%%%%%%%%%%%%%%%%
%
%%%%%%%%%%%%%%%%%
%
%%%%%%%%%%%%%%%%%%%%%%%%%%%%%%%%%%%%%%%%%%% Acknowledgement %%%%%
%
%%%%%%%%%%%%%%%%%%%%%%%%%%%%%%%%%%%%%%%%%%%%%%%%%%%%%%%%%%%%%%%%%
\section*{Acknowledgement}
%%%%%%%%%%%%%%%%%%%%%%%%%%%%%%%%%%%%%%%%%%%%%%%%%%%%%%%%%%%%%%%%%
%
%%%%%%%%%%%%%%%%%%%%%%%%%%%%%%%%%%%
The author is grateful to Professor Chuichiro Hayashi for his 
helpful comments and encouragements.
%%%%%%%%%%%%%%%%%%%%%%%%%%%%%%%%%%%
%
%%%%%%%%%%%%%%%%
\par\vspace{5mm}
%%%%%%%%%%%%%%%%
%
%
%%%%%%%%%%%%%%%%%%%%%%%%%%%%%%%%%%%%%%%%%%%%%%%%%%%%%%%%%%%%%%%%

%%%%%%         
%
%\addtocounter{figure}{-27}\Alph{figure}
%
%%%%%%
\begin{figure}[htbp]\begin{center}
\includegraphics[trim=0mm 0mm 0mm 0mm, width=.8\linewidth]
{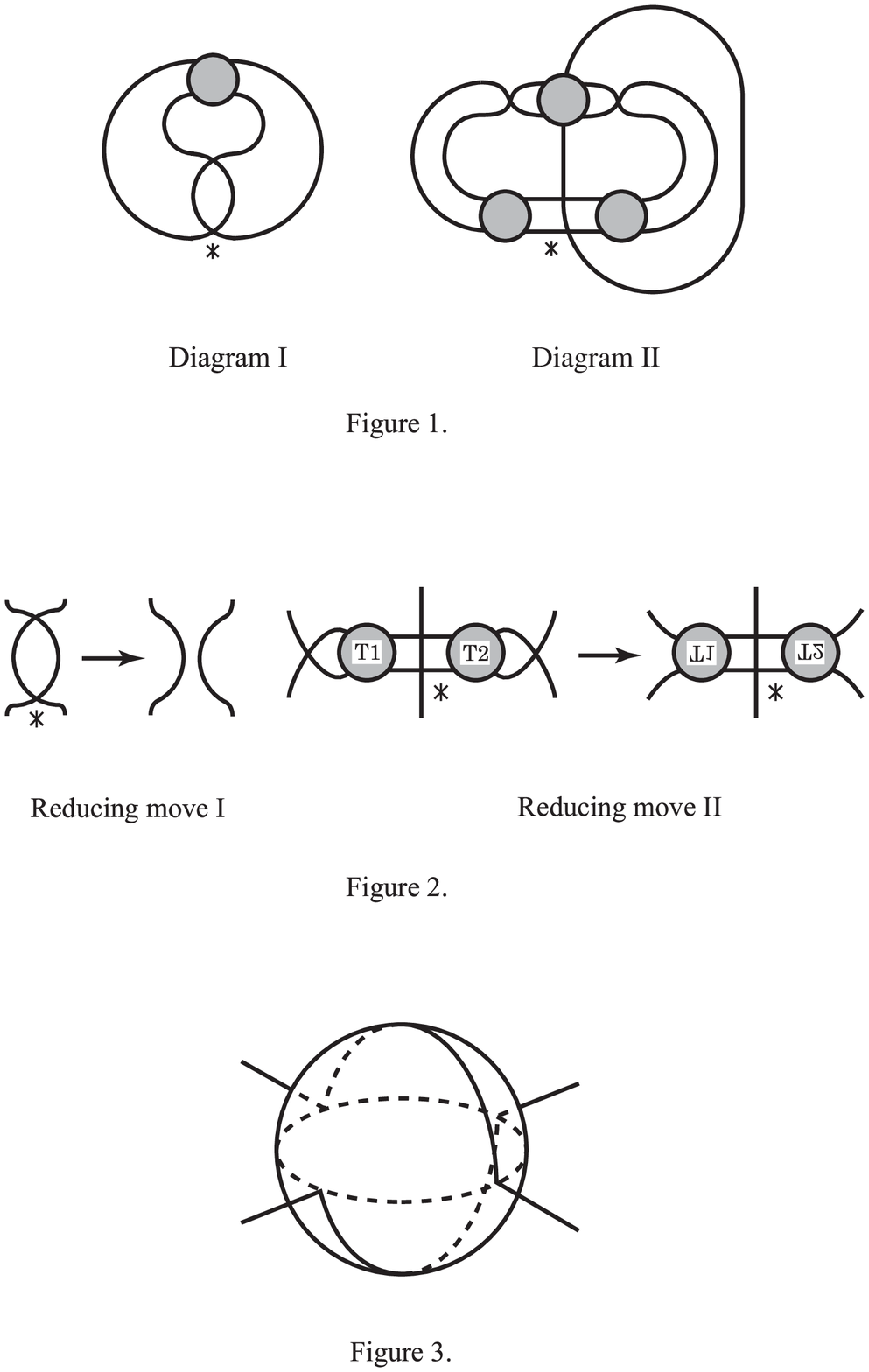}\end{center}\end{figure}
%%%%%%
\begin{figure}[htbp]\begin{center}
\includegraphics[trim=0mm 0mm 0mm 0mm, width=.8\linewidth]
{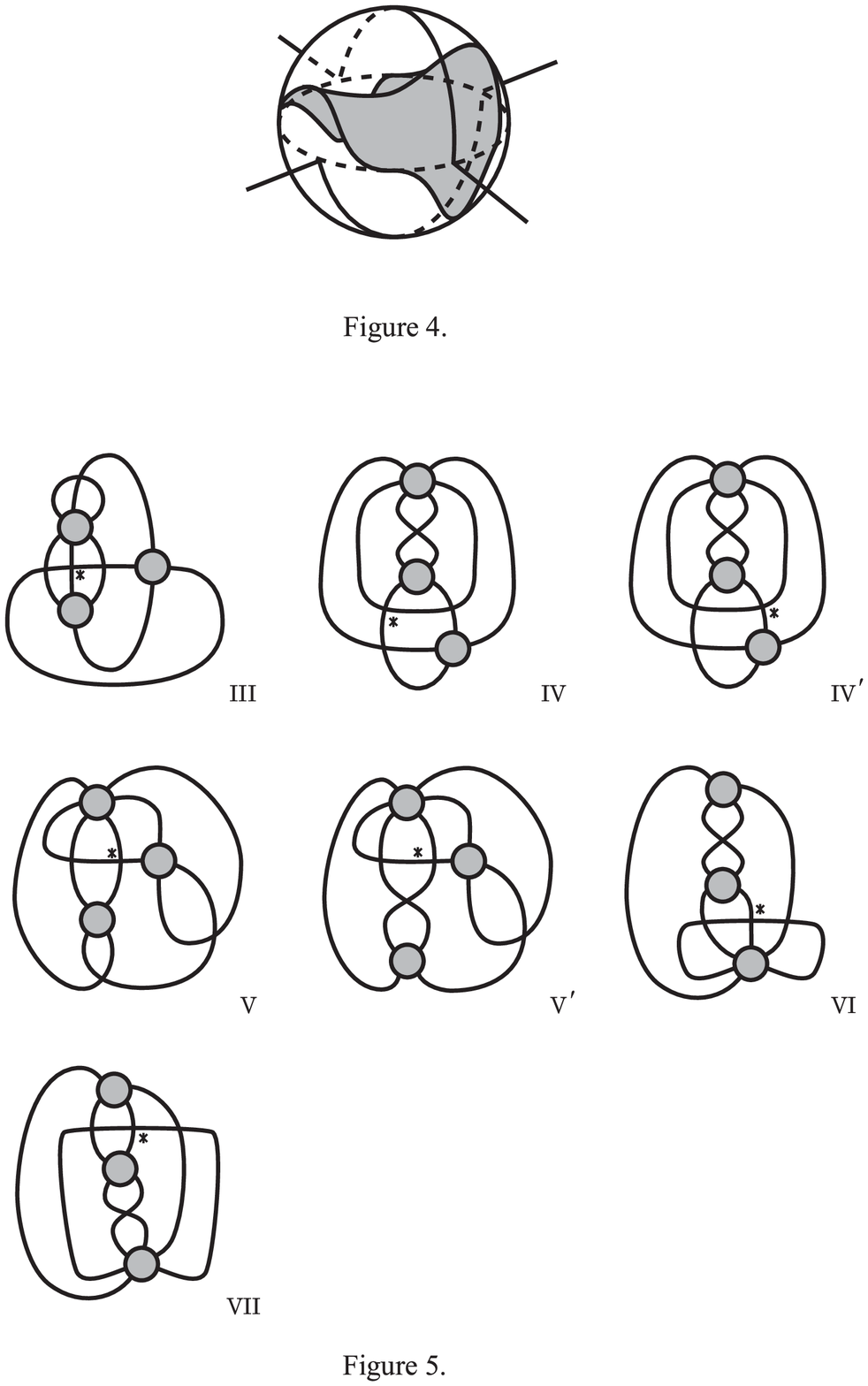}\end{center}\end{figure}
%%%%%%
\begin{figure}[htbp]\begin{center}
\includegraphics[trim=0mm 0mm 0mm 0mm, width=.7\linewidth]
{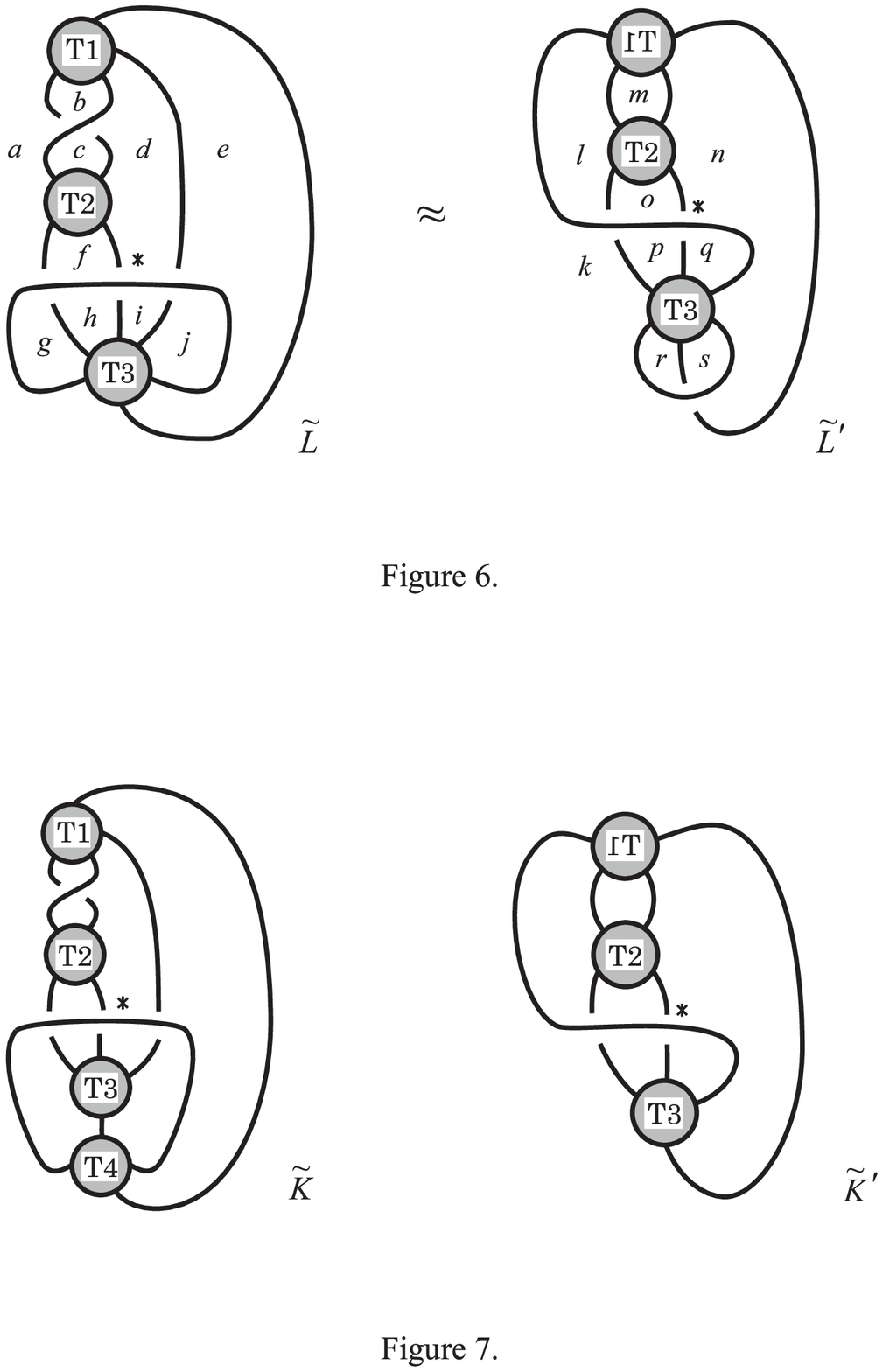}\end{center}\end{figure}
%%%%%%
\begin{figure}[htbp]\begin{center}
\includegraphics[trim=0mm 0mm 0mm 0mm, width=.6\linewidth]
{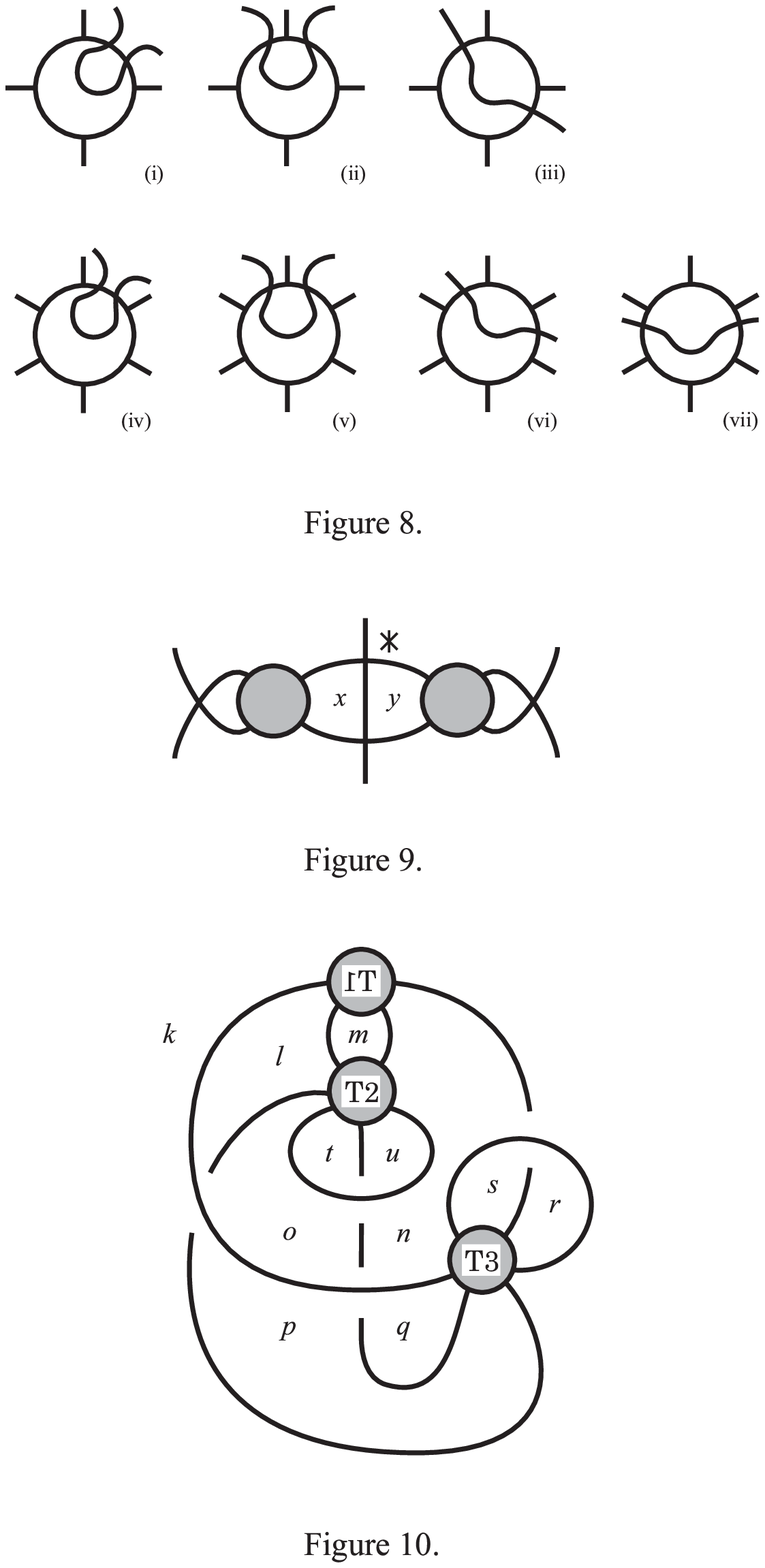}\end{center}\end{figure}
%%%%%%
\begin{figure}[htbp]\begin{center}
\includegraphics[trim=0mm 0mm 0mm 0mm, width=.7\linewidth]
{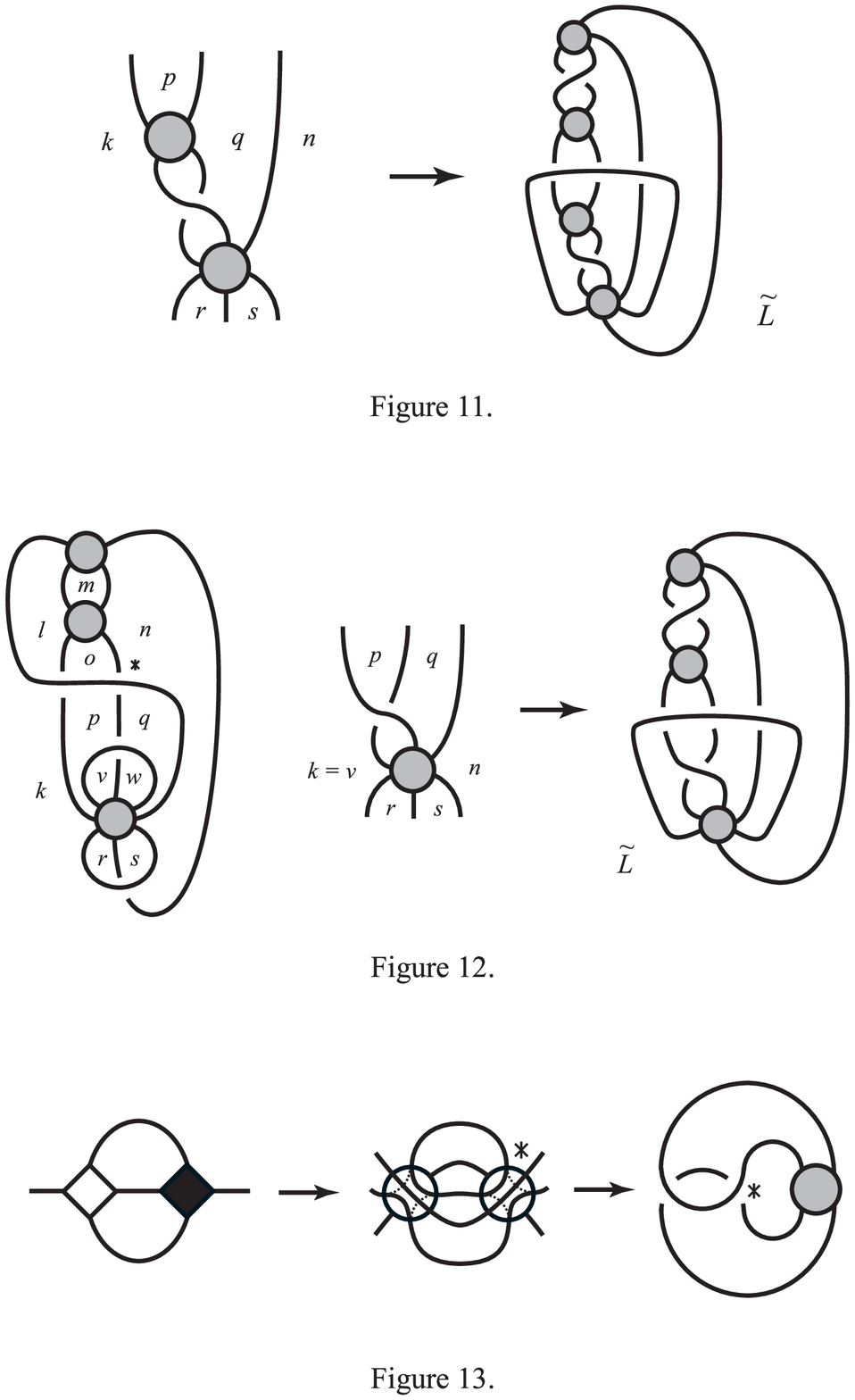}\end{center}\end{figure}
%%%%%%
\begin{figure}[htbp]\begin{center}
\includegraphics[trim=0mm 0mm 0mm 0mm, width=.8\linewidth]
{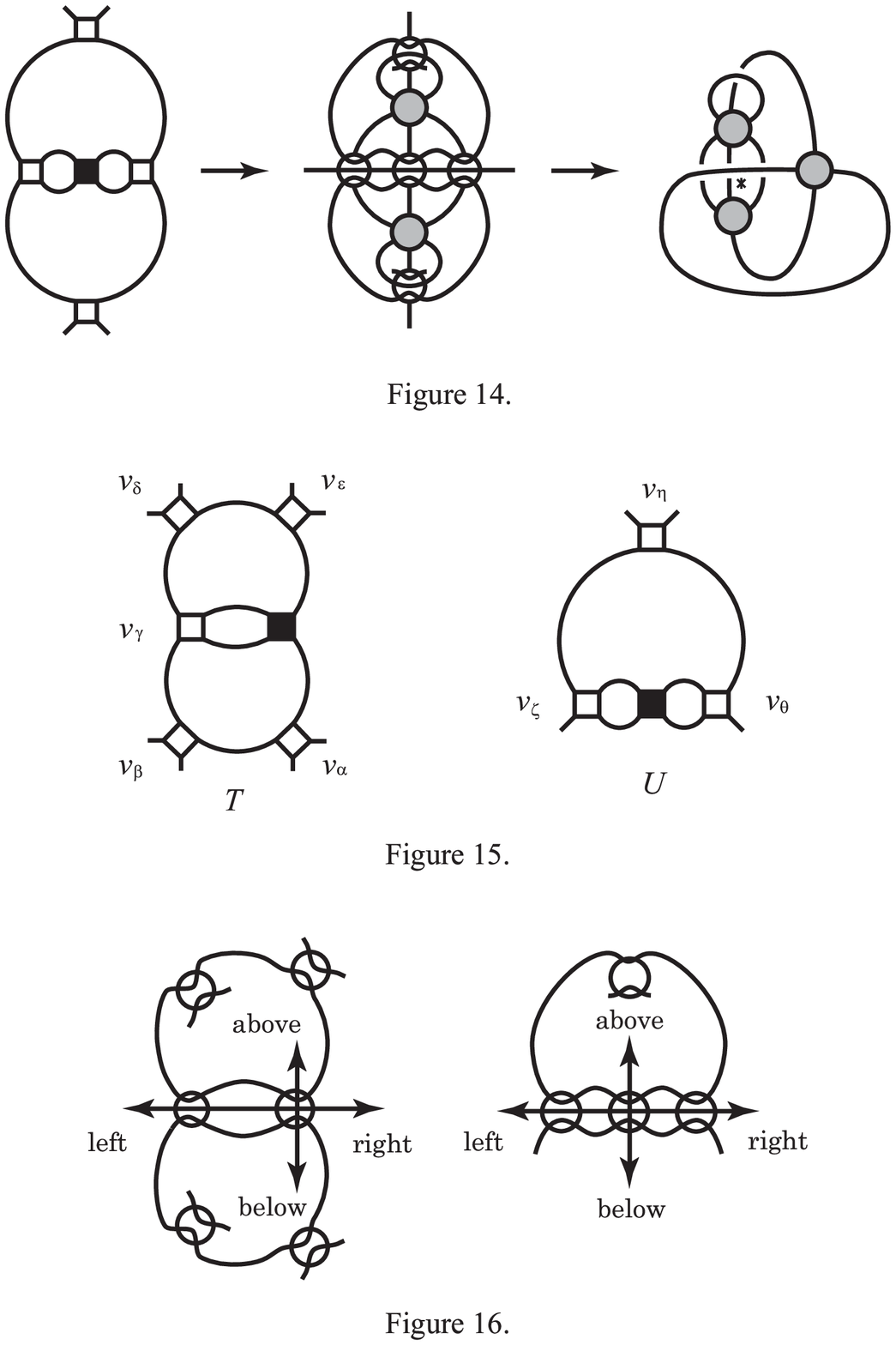}\end{center}\end{figure}
%%%%%%
\begin{figure}[htbp]\begin{center}
\includegraphics[trim=0mm 0mm 0mm 0mm, width=.8\linewidth]
{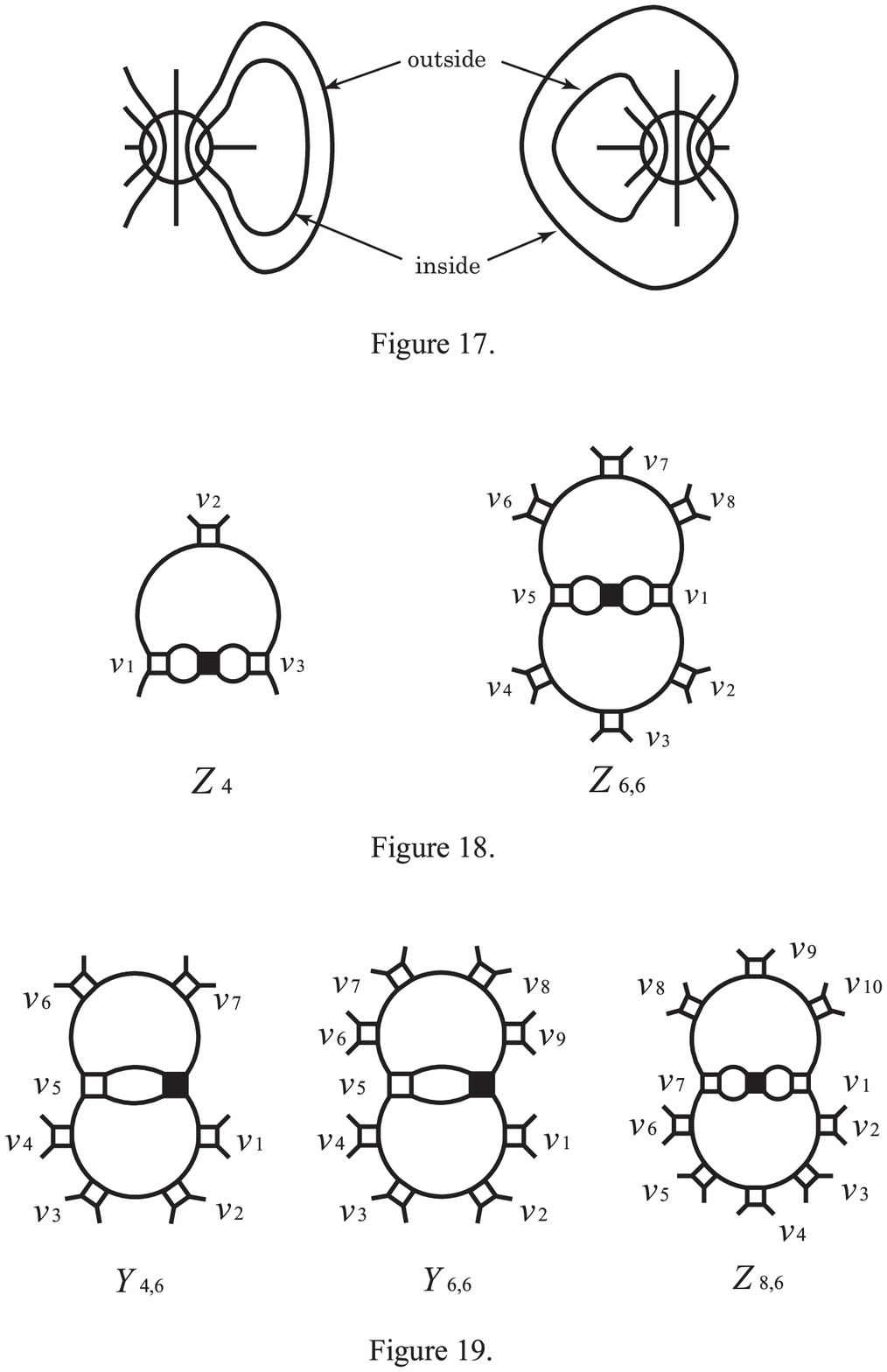}\end{center}\end{figure}
%%%%%%
\begin{figure}[htbp]\begin{center}
\includegraphics[trim=0mm 0mm 0mm 0mm, width=.8\linewidth]
{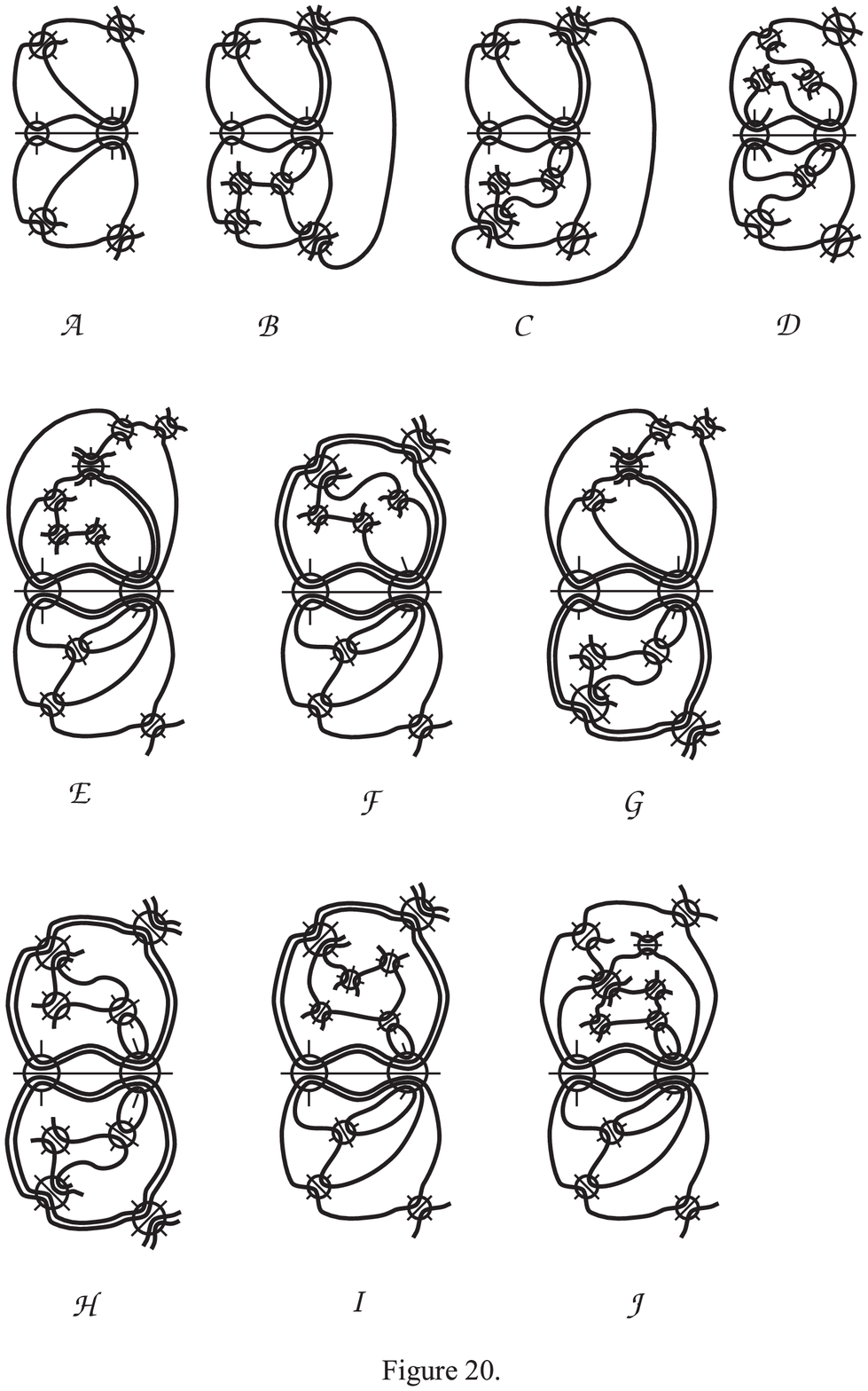}\end{center}\end{figure}
%%%%%%
\begin{figure}[htbp]\begin{center}
\includegraphics[trim=0mm 0mm 0mm 0mm, width=.8\linewidth]
{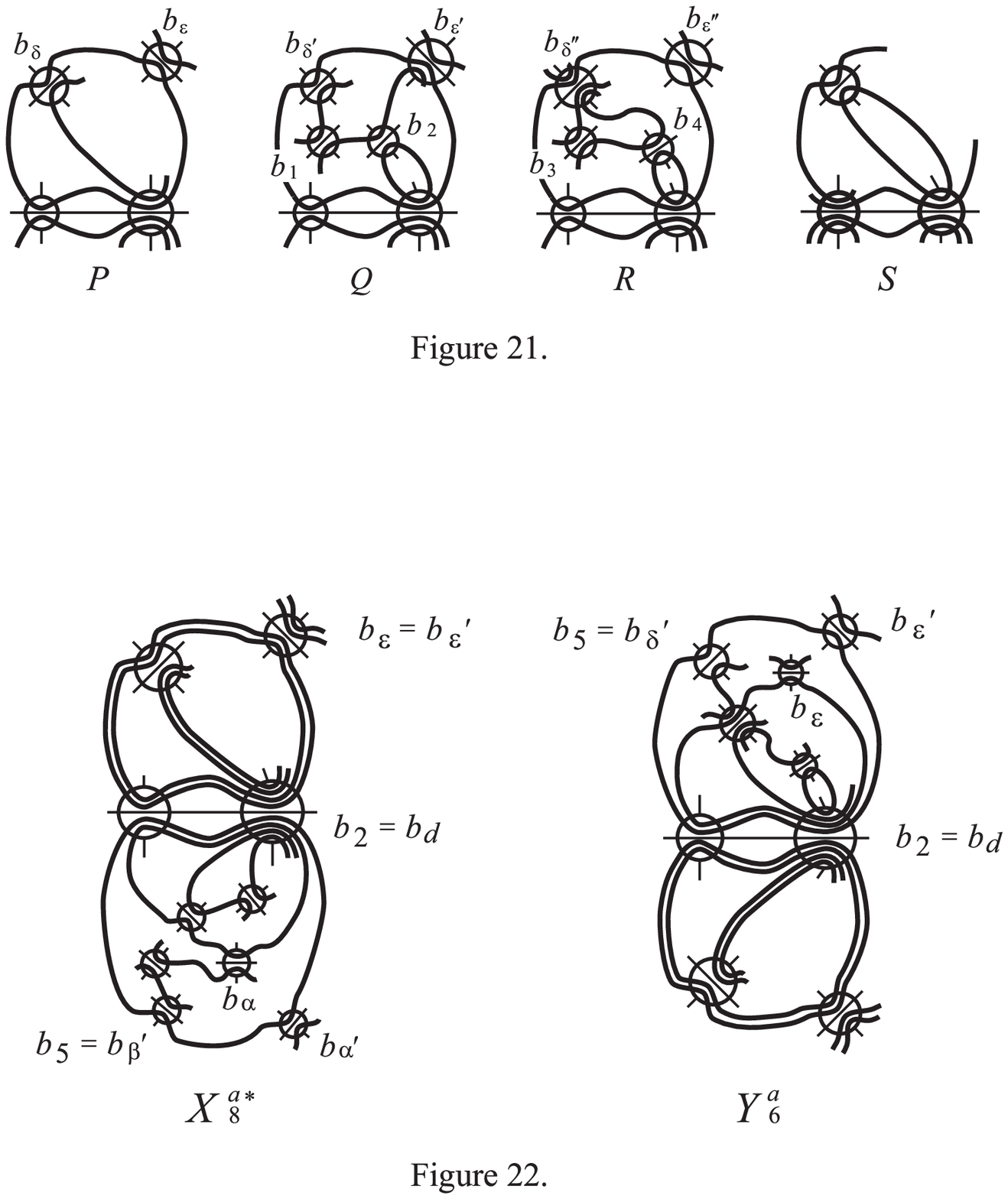}\end{center}\end{figure}
%%%%%%
\begin{figure}[htbp]\begin{center}
\includegraphics[trim=0mm 0mm 0mm 0mm, width=.8\linewidth]
{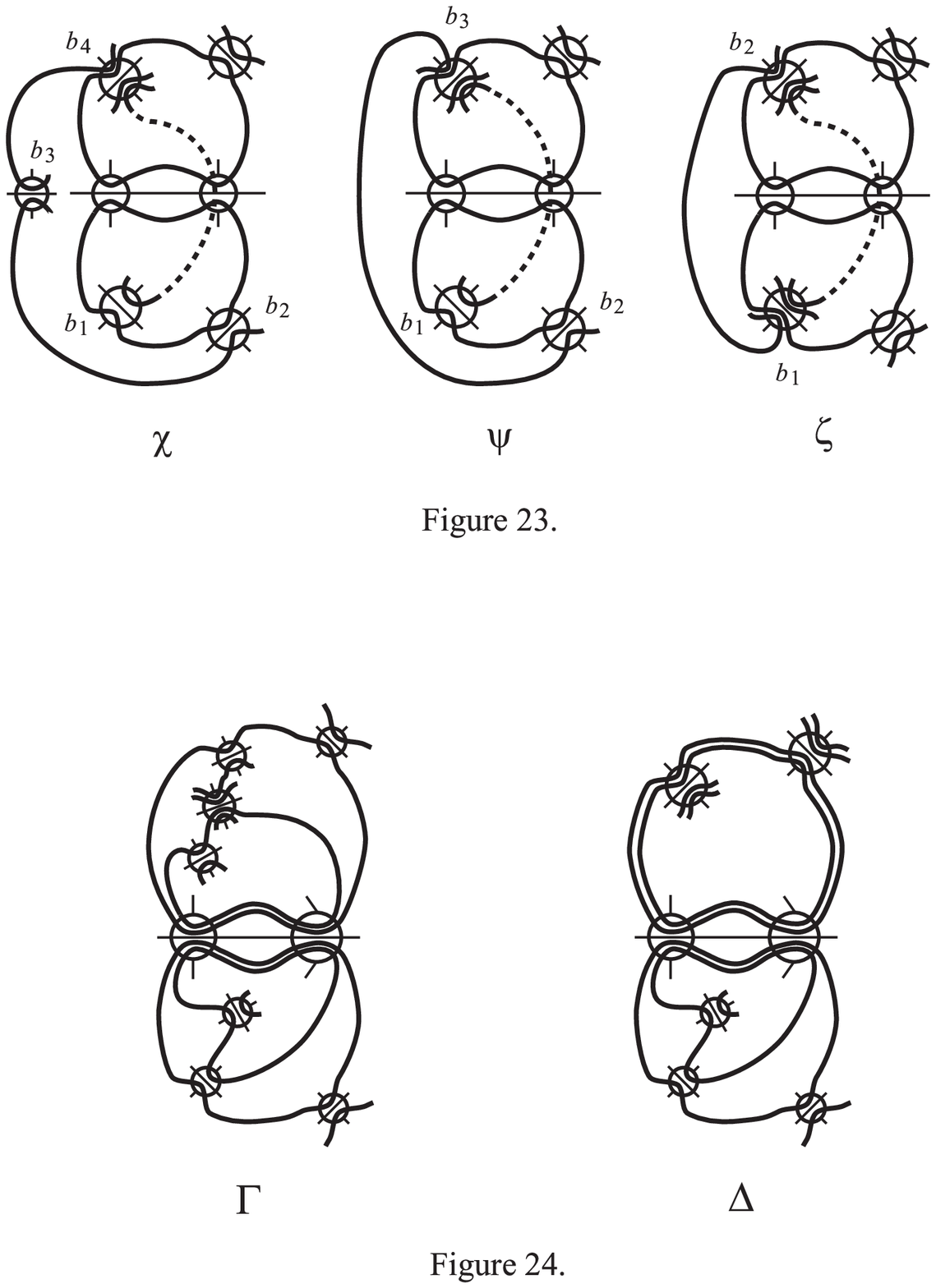}\end{center}\end{figure}
%%%%%%
\begin{figure}[htbp]\begin{center}
\includegraphics[trim=0mm 0mm 0mm 0mm, width=.8\linewidth]
{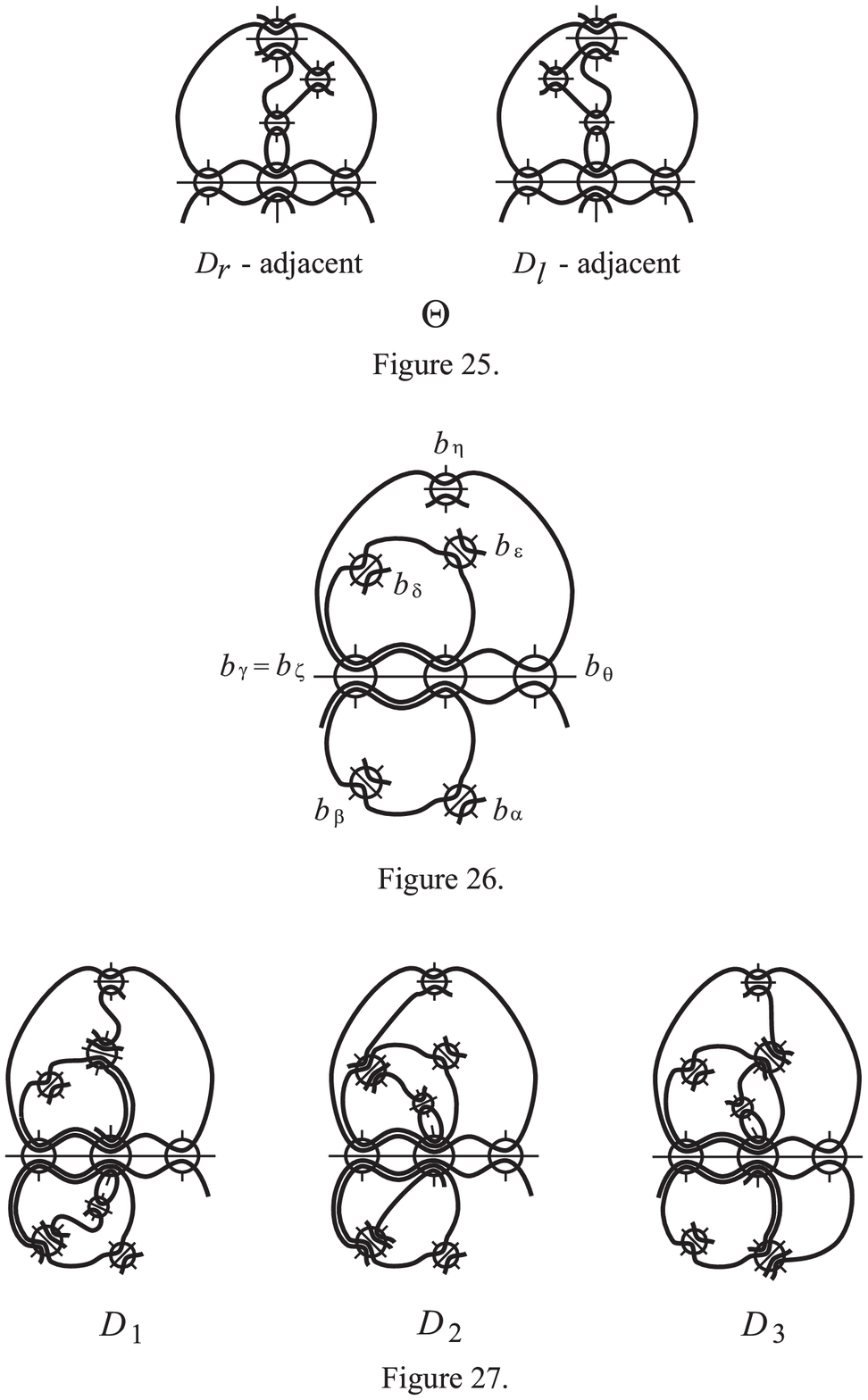}\end{center}\end{figure}
%%%%%%
\end{document}